\documentclass[final]{siamltex}
\usepackage{graphicx,color}
\usepackage{diagbox}
\usepackage{multirow}
\usepackage{subfig}
\usepackage{epstopdf}
\usepackage{comment}
\usepackage{braket}
\usepackage{caption}
\usepackage{adjustbox}
\usepackage{amsmath, amssymb, amsfonts}
\usepackage{subfig}
\usepackage{xcolor}
\usepackage{bbm}
\usepackage{enumerate}
\usepackage{bm}

\usepackage{algorithm}
\usepackage{algpseudocode}

\algrenewcommand\textproc{}

\newcommand\pfrac[2]{\frac{\partial #1}{\partial #2}}
\newcommand\bbRm[2]{\mathbb{R}^{#1 \times #2}}
\newcommand\figref[1]{Fig.\ \ref{#1}}
\newcommand\secref[1]{Section \ref{#1}}

\newcommand\appref[1]{Appendix \ref{#1}}

\newcommand{\bvec}[1]{\bm{#1}}

\newcommand{\vr}{\bvec{r}}

\newcommand{\vi}{\bvec{i}}
\newcommand{\vj}{\bvec{j}}

\newcommand{\vx}{\bvec{x}}

\newcommand{\abs}[1]{\lvert#1\rvert}

\DeclareMathOperator*{\argmin}{arg\,min}

\newcommand{\bbP}{\mathbb{P}}
\newcommand{\bbR}{\mathbb{R}}

\newcommand{\ccA}{\mathcal{A}}
\newcommand{\ccB}{\mathcal{B}}
\newcommand{\ccC}{\mathcal{C}}
\newcommand{\ccD}{\mathcal{D}}
\newcommand{\ccG}{\mathcal{G}}

\newcommand{\ccO}{\mathcal{O}}
\newcommand{\ccP}{\mathcal{P}}
\newcommand{\ccQ}{\mathcal{Q}}

\global\long\def\R{\mathbb{R}}

\global\long\def\Rd{\mathbb{R}^{d}}
\global\long\def\E{\mathbb{E}}
\global\long\def\P{\mathbb{P}}

\global\long\def\ra{\rightarrow}

\numberwithin{equation}{section}
\numberwithin{figure}{section}
\newtheorem{thm}{\protect\theoremname}

\newtheorem{defn}[thm]{Definition}

\providecommand{\corollaryname}{Corollary}
\providecommand{\lemmaname}{Lemma}
\providecommand{\propositionname}{Proposition}
\providecommand{\remarkname}{Remark}
\providecommand{\theoremname}{Theorem}

\newcommand{\OAB}{\Omega\backslash(A\cup B)}

\allowdisplaybreaks

\title{Committor functions via tensor networks}

\author{Yian Chen\thanks{Department of Statistics, University of Chicago, Illinois, IL 60637, USA. Email: {\tt yianc@uchicago.edu}}
\and
Jeremy Hoskins\thanks{Department of Statistics, University of Chicago, Illinois, IL 60637, USA. Email: {\tt jeremyhoskins@uchicago.edu}}
\and
Yuehaw Khoo\thanks{Department of Statistics, University of Chicago, Illinois, IL 60637, USA. Email: {\tt ykhoo@uchicago.edu}}
\and
Michael Lindsey\thanks{Department of Mathematics, Courant Institute of Mathematical Sciences, New York University, New York, NY 10012, USA.  Email: \texttt{michael.lindsey@cims.nyu.edu}}
}

\begin{document}

\maketitle

\begin{abstract}
    We propose a novel approach for computing committor functions, which describe transitions of a stochastic process between metastable states. The committor function satisfies a backward Kolmogorov equation, and in typical high-dimensional settings of interest, it is intractable to compute and store the solution with traditional numerical methods. By parametrizing the committor function in a matrix product state/tensor train format and using a similar representation for the equilibrium probability density, we solve the variational formulation of the backward Kolmogorov equation with linear time and memory complexity in the number of dimensions. This approach bypasses the need for sampling the equilibrium distribution, which can be difficult when the distribution has multiple modes. Numerical results demonstrate the effectiveness of the proposed method for high-dimensional problems.
\end{abstract}



\section{Introduction}
Understanding rare transitions between metastable states of a high-dimensional stochastic processes is a problem of great importance in the applied sciences. Examples of interesting transition events include chemical reactions, nucleation events during phase transitions, and conformational changes of molecules \cite{chem-tran-1,chem-tran-2,nuc-1,nuc-2,conf-change-1,conf-change-2}. In such complex systems, the dynamics linger near metastable states for long waiting periods, punctuated by sudden jumps from one metastable state to another. One important tool for describing transition events is transition path theory \cite{transition-1,transition-2,transition-3, transition-info-1}, where the  committor function plays a central role. The committor function measures the probability that the process hits a certain metastable state of the system before another and can 
be viewed as the solution of a backward Kolmogorov equation.

Computing the committor function in high-dimensional settings is a formidable task. Traditional numerical methods such as finite difference and finite element methods become prohibitively expensive in even moderate dimensions. To overcome the curse of dimensionality, significant efforts have been expended to apply deep learning framework to solve for high dimensional partial differential equations \cite{deepritz,deep-pde-1,deep-pde-2,deep-pde-3,pinn}. Most recently \cite{NN-committor-1,LiLinRen2019,NN-committor-2} have suggested representing the committor function using neural networks. Some of these 
approaches rely on sampling from the equilibrium distribution and so work well when the transitions are easily observed. For, e.g., chemical systems at low temperature the committor function can change sharply between the two metastable states, and transitions are rare and difficult to sample. To address this problem, \cite{Rostkoff2021} has proposed an adaptive importance sampling scheme. Meanwhile, the dynamical Galerkin framework for computing committor functions \cite{dga}, which represents functions in a basis rather than as neural networks, approaches the sampling problem by initializing short trajectories uniformly according to known reaction coordinates.

Tensor network methods \cite{tt-quantum-1,tensor-diagram-intro,tt-quantum-3,tt-quantum-4,tt-quantum-5, tensorGalerkin2011} have emerged as an alternative to neural networks as a tool for high-dimensional problems in modern quantum physics and beyond. Typical tensor decomposition methods include tensor trains \cite{tt-decomposition} (also known as linear tensor networks or matrix product states \cite{mps-1,mps-2,mps-3}), the CP decomposition \cite{cp-decomposition}, and the Tucker/Hierarchical Tucker decomposition \cite{cp-decomposition,tucker-decomposition,h-tucker}. These methods approximate tensors in compressed, structured formats that enable efficient linear algebra operations. More details can be found in \cite{tt-review-1,tt-review-2,tt-review-3,tt-review-4}. Moreover, tensor network methods have also been applied to solve for high-dimensional partial differential equations \cite{tensor-pde-1,tensor-pde-2}. 

In this paper we propose applying the tensor method to compute the committor functions based on matrix product states/tensor trains. Specifically, we approximate both the equilibrium probability distribution and the committor function using tensor trains, achieving good performance even for high-dimensional problems in the low-temperature regime.

The rest of the paper is organized as follows. The committor function and its properties are reviewed in \secref{sec:background}. Therein we also explain how the 
boundary condition can be accommodated within our 
tensor format and provide a summary of
relevant tensor network methods. We introduce the key ingredients of our proposed method in \secref{sec:method}. Numerical experiments for two representative classes of examples are presented in \secref{sec:numerical-tests}, demonstrating the accuracy and efficiency of the proposed algorithm. Finally, in \secref{sec:conclusion} we summarize our findings.

\section{Background and preliminaries}
\label{sec:background}

In this section we first review the motivation for computing committor functions, and summarize challenges and recent advances relevant to this task. We then briefly discuss tensor train decompositions, introduce the basic tensor operations, and define relevant associated notations used in this work. Throughout the paper, we use MATLAB notation for multidimensional array indexing. 

\subsection{Committor functions}

The underyling stochastic process of interest is the overdamped Langevin process, defined by 
\begin{align}
    d \bm{X}_t = -\nabla V(\bm{X}_t) \, dt + \sqrt{2\beta^{-1}} \, d\bm{W}_t,
    \label{eq:langevin}
\end{align}
where $\bm{X}_t \in \Omega\subset\Rd$ is the state of the system, $V:\Omega\subset\Rd\rightarrow\bbR$ is a smooth potential energy function, $\beta=1/T$ is the inverse of the temperature $T$, and $\bm{W}_t$ is a  $d$-dimensional Wiener processs. If the potential energy function $V$ is confining for $\Omega$ (see, e.g., \cite[Definition 4.2]{bhattacharya2009stochastic}), then one can show that the equilibrium probability distribution of the Langevin dynamics \eqref{eq:langevin} is the Boltzmann-Gibbs distribution
\begin{align}
    p(\vx) = \frac{1}{Z_\beta}\exp(-\beta V(\vx)).
    \label{eq:equi-measure}
\end{align}
where $Z_\beta=\int_{\Omega} \exp(-\beta V(\vx))\, d\vx$ is the partition function. We are interested in the transition between two simply connected domains $A, B \subset \Omega$ with smooth boundaries. The associated committor function $q:\Omega \to [0,1]$ is defined by
\begin{align}
    q(\vx) = \bbP(\tau_B<\tau_A\ |\ \bm{X}_0=\vx),
\end{align}
where $\tau_A$ and $\tau_B$ are the hitting times for the sets $A$ and $B$, respectively. The committor function $q$ provides a useful statistical description of properties such as the density and probability of reaction trajectories \cite{transition-3,transition-info-1,transition-info-2}. However, computing the committor function can be a
formidable task since it involves solving the following (possibly high-dimensional) backward Kolmogorov equation with Dirichlet boundary conditions:
\begin{align}
\label{eq:kolm}
    -\beta^{-1}\Delta q(\vx) + \nabla V(\vx) \cdot \nabla q(\vx) = 0 \ \text{in}\ \OAB,\ q(\vx)|_{\partial A} = 0,\ q(\vx)|_{\partial B} = 1.
\end{align}
For high-dimensional problems, traditional methods such as finite difference and finite element discretization are intractable. Numerous alternative methods can effectively approximate the committor function under the assumption that the transition paths from $A$ to $B$ are localized in a quasi-one-dimensional reaction tube or low-dimensional manifold. For example, the finite temperature string method \cite{finite-temp-string,revisit-finite-temp-string} approximates the isosurfaces of the committor function with hyperplanes normal to the most probable transition paths. The diffusion map approach \cite{diffusion-map-committor} aims to obtain the committor function on a set of points by applying point cloud discretization to the generator $L = -\beta^{-1}\Delta + \nabla V(\vx) \cdot \nabla$. The method presented in \cite{point-committor} improves the diffusion map approach by discretizing $L$ using a finite element method on local tangent planes of the point cloud.

Recently, \cite{NN-committor-1,LiLinRen2019,NN-committor-2,Rostkoff2021} obtained the committor function by way of the variational problem
\begin{align}
    \argmin_q \int_{\Omega} \abs{\nabla q(\vx)}^2 p(\vx) \, d\vx,\ \ q(\vx)|_{\partial A} = 0,\ q(\vx)|_{\partial B} = 1,
    \label{eq:cons-varational}
\end{align}
for which \eqref{eq:kolm} is the Euler-Lagrange equation, by 
optimizing a a neural network parametrization for the committor function $q$. In order to obtain an unconstrained optimization problem, the boundary conditions are enforced in \cite{NN-committor-1,NN-committor-2,Rostkoff2021} by adding two extra penalty terms, as in 
\begin{align}
    \argmin_{q} \int_{\Omega} \abs{\nabla q(\vx)}^2 p(\vx) \, d\vx + \rho \int_{\partial A} q(\vx)^2 p_{\partial A}(\vx) \, d\vx + \rho \int_{\partial B} (q(\vx) - 1)^2 p_{\partial B}(\vx) \, d\vx,
    \label{eq:uncons-varational}
\end{align}
where $p_{\partial A}$ and $p_{\partial B}$ define probability measures supported on the boundaries $\partial A$ and $\partial B$ respectively. In all of these works, the objective is evaluated and optimized via stochastic sampling. By contrast, in this work, we propose to represent the committor function $q$ in a tensor train format, which will 
allow for optimization via stable and efficient deterministic linear algebra 
operations.

Since the potential function $V$ is confining, we can effectively restrict our domain of interest to a bounded subset of $\Omega$. Outside of this subset, the density is small and contributes only negligibly contributes to the variational cost \eqref{eq:cons-varational}. For simplicity we shall identify $\Omega$ with this subset and assume $\Omega=\Omega_1\times\Omega_2\times\dots\Omega_d$ where each $\Omega_i\in\bbR$ is a bounded subset. 

\subsection{Soft boundary condition}
Unfortunately, the formulation \eqref{eq:uncons-varational} is
not immediately amenable to optimization within a tensor format for $q$. The reason is that 
the surface measures on $\partial A$ and $\partial B$ cannot themselves be identified with functions on $\Omega$, much 
less as functions that can be compressed in tensor format, so 
the penalty terms cannot simply be viewed as inner products of tensor trains.

Therefore we instead consider an objective of the form 
\begin{align}
    \argmin_{q} \int_{\Omega} \abs{\nabla q(\vx)}^2 p(\vx) \, d\vx + \rho \int_{\Omega} q(\vx)^2 p_{A}(\vx) \, d\vx + \rho \int_{\Omega} (q(\vx) - 1)^2 p_{B}(\vx) \, d\vx.
    \label{eq:soft-varational}
\end{align}
Here $p_A$ and $p_B$ are probability densities that are absolutely continuous with respect to the Lebesgue measure on $\Omega$.

In fact, we show in Appendix~\ref{app:soft} that the exact optimizer of \eqref{eq:soft-varational} admits a probabilistic interpretation similar to that of the usual committor function. As such we call the optimizer a `soft committor function.' Specifically, the interpretation is based on a modification of the Langevin dynamics~\eqref{eq:langevin} in which one augments the state 
space $\Omega$ with two `cemetery states' $c_A$ and $c_B$. The process jumps randomly to these two states with instantaneous jump rates $\frac{\rho \,\cdot\, p_A}{\beta \, \cdot\, p}$ and
$\frac{\rho\,\cdot \, p_B}{\beta \,\cdot \,  p}$, respectively. The soft committor function evaluates the probability that the modified process hits $c_B$ before $c_A$.

Evidently, when $\rho$ is large and $p_A$ and $p_B$ concentrate near $\partial A$ and $\partial B$, respectively, the soft committor function can be viewed as an approximation 
of the ordinary committor function. In fact, if $p_A$ and $p_B$ are Gaussian densities, then in high dimensions \cite{foundation-of-ds}, $p_A$ and $p_B$ each weakly approximates a uniform measure on a suitable hypersphere. 
This is convenient because $A$ and $B$ are often chosen to be balls. More details on the construction of $p_A$ and $p_B$ will be provided below in Section \ref{sec:method}.

For simplicity, in what follows we will simply refer to soft committor functions as committor functions.

\subsection{Tensors and tensor networks}

In this subsection we summarize the basic tensor operations used in this work. In particular, for ease of exposition we introduce tensor network diagram notation, which provides a convenient way of visually describing tensor operations. We also introduce the matrix product state/tensor train format for parametrizing high-dimensional functions.

\begin{figure}[!htb]
    \centering
    \subfloat[]{{\includegraphics[width=0.40\textwidth]{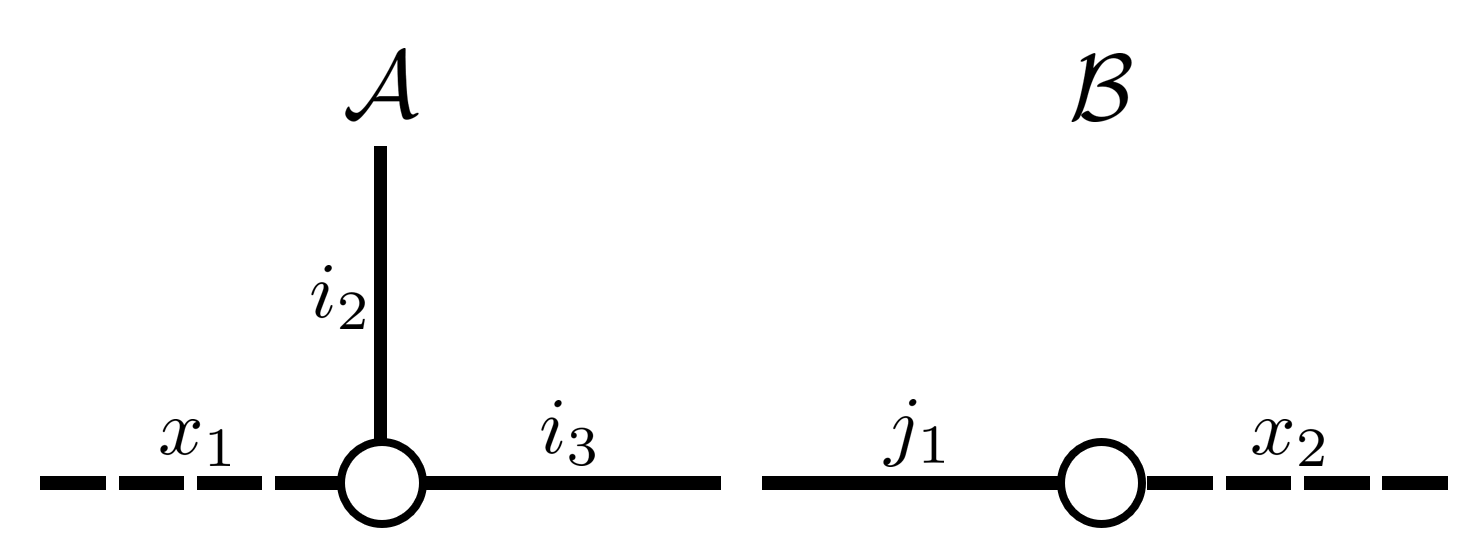}}}\quad \quad \quad \quad
    \subfloat[]{{\includegraphics[width=0.20\textwidth]{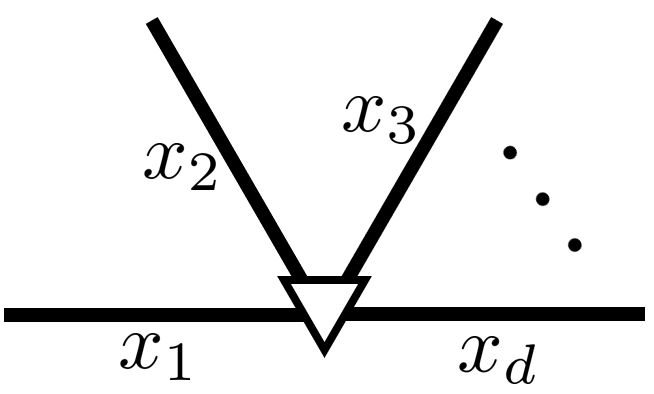}}}
    \caption{(a) Tensor diagrams for a $3$-tensor $\ccA$ and a $2$-tensor $\ccB$. (b) Tensor diagram for a $d$-dimensional Kronecker Delta node. Solid lines correspond to discrete variables, and dashed lines correspond to continuous variables.} 
    \label{fig:tensor-and-delta}
\end{figure}

In tensor diagrams, a tensor is represented by a node, where the number of incoming legs indicates the dimensionality of the tensor, i.e., the number of indices/arguments. There are two types of leg: legs indicating continuous arguments are denoted by dashed lines, and legs indicating discrete indices are denoted by solid lines. For example, \figref{fig:tensor-and-delta} (a) shows the tensor diagram for a $3$-tensor $\ccA$ and a $2$-tensor $\ccB$, which can be viewed as two functions
\begin{equation}
    \ccA(x_1,i_2,i_3),\quad \ccB(j_1,x_2),
\end{equation}
respectively, where $x_1,x_2$ are continuous variables and $i_2,i_3,j_1$ are discrete variables.

We also define the multi-dimensional Kronecker delta tensor (depicted by an inverted triangular node as in \figref{fig:tensor-and-delta} (b)):
\begin{align}
    \delta(x_1,x_2,\dots,x_d) = 
    \begin{cases}
        1 & \text{if } x_1=x_2=\dots=x_d.\\
        0 & \text{otherwise}.
    \end{cases}
\end{align}
By a slight abuse of notation, we will use the same symbol to represent the appropriate Dirac delta function when the legs represent continuous variables.

\begin{figure}[!htb]
    \centering
    \subfloat[]{{\includegraphics[height=0.125\textwidth]{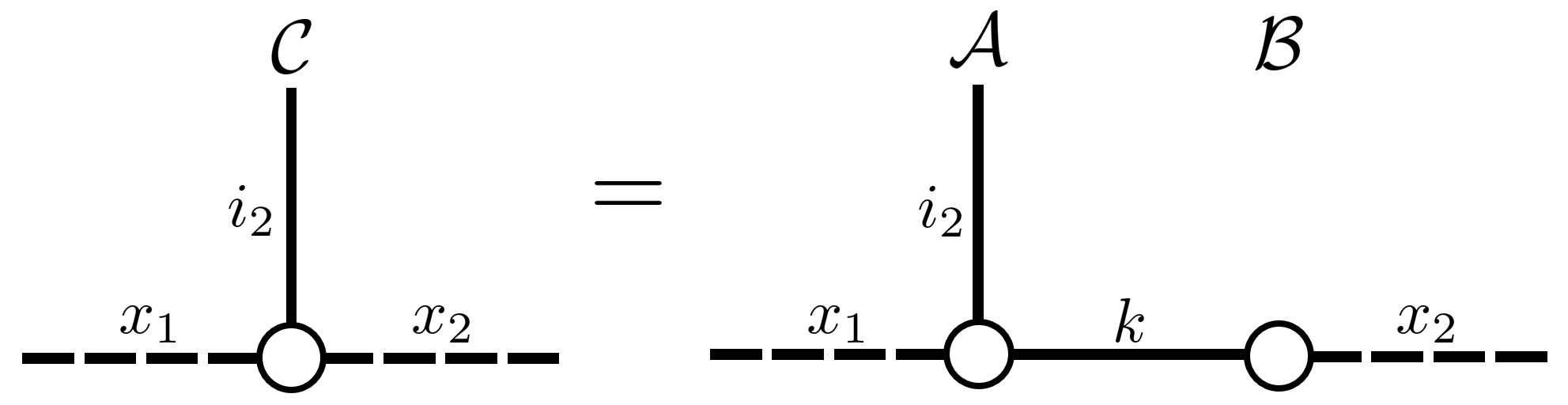}}}\\ \vspace{2 mm}
    \subfloat[]{{\includegraphics[height=0.125\textwidth]{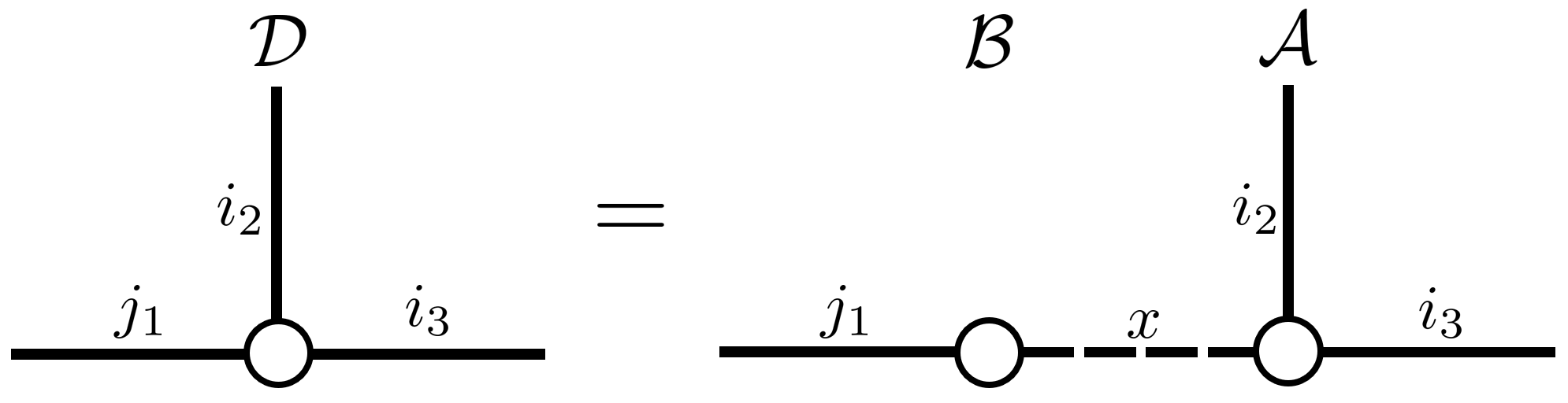}}}
    \caption{(a) Tensor contraction of discrete legs. (b) Tensor contraction of continuous legs.} 
    \label{fig:tensor-contraction}
\end{figure}

Next we describe a key operation called tensor contraction. This operation is indicated visually by joining legs from different tensors. For example, in \figref{fig:tensor-contraction} (a), the third leg of $\ccA$ is joined with the first leg of $\ccB$. This corresponds to the computation
\begin{align}
    \ccC(x_1,i_2,x_2) &= \sum_{k} \ccA(x_1,i_2,k) \ccB(k,x_2),
    \label{eq:tensor-contraction-1}
\end{align}
where it is implicitly assumed that the indices of the joined legs have the same range. Here $x_1,x_2$ are continuous variables, and $i_2$ is a discrete variable.

Tensor contraction can be defined for continuous legs as well. For example, in \figref{fig:tensor-contraction} (b), continuous  legs of $\ccA$ and $\ccB$ are contracted, corresponding to the operation
\begin{align}
    \ccD(j_1,i_2,i_3)
    &= \int_{\Omega_0} \ccB(j_1,x) \ccA(x,i_2,i_3) \, dx,
    \label{eq:tensor-contraction-2}
\end{align}
for some suitable domain $\Omega_0$, which is implicitly assumed to be the domain of both joined legs. The resulting tensor $\ccD$ is a $3$-tensor with only discrete legs.

A tensor network diagram consists of a collection of individual tensor diagrams with some pairs of legs joined, i.e., contracted. The contracted legs correspond to the so-called `internal indices' for the tensor network, while the uncontracted legs correspond to `external indices,' which are the indices remaining after all of the indicated contractions have been performed.

Next we introduce several low-complexity tensor networks and their corresponding diagrams. A matrix product state (MPS) or tensor train (TT) is a factorization of a $d$-tensor into a chain-like product of $3$-tensors. Such a factorization allows one to approximate high-dimensional tensors and manipulate them efficiently, typically with $\ccO(d)$ time and memory complexity.
~\\
\begin{defn}
    \label{def:mps}
    Let $\ccA\in\bbR^{n_1\times n_2\times\dots\times n_d}$ be a $d$-tensor, with entries indexed by $(i_1,i_2,\dots,i_d)$. Then we say that $\ccA$ is a MPS/TT with ranks $\vr=(r_0,\dots,r_{d})$, where we fix $r_0 = r_d = 1$ by convention, if one can write
    \begin{align}
        \ccA(i_1,i_2,\dots,i_d)\  =& \ \  \sum_{\alpha_0 = 1}^{r_0} \cdots \sum_{\alpha_{d} = 1}^{r_{d}} \ccG_1(\alpha_0, i_1,\alpha_1)\ccG_2(\alpha_1, i_2, \alpha_2)\dots\ccG_d(\alpha_{d-1}, i_d, \alpha_d)\nonumber\\
        =&\ \  \ccG_1(:, i_1,:)\ccG_2(:, i_2, :)\cdots \ccG_d(:, i_d, :)
        \label{eq:TT}
    \end{align}
    for all $(i_1,i_2,\dots,i_d)$.
    Here $\ccG_k(:, i_k, :)\in\bbRm{r_{k-1}}{r_k}$ is viewed as a matrix for each $k = 1,\ldots, d$, and the matrix product in \eqref{eq:TT} is a $1\times 1$ matrix, i.e., a scalar value. The 3-tensor $\ccG_k \in \R^{r_{k-1} \times n_k \times r_{k}}$ is called the $k$-th tensor core of $\ccA$.
\end{defn}
~\\

\begin{figure}[!htb]
    \centering
    \subfloat[]{{\includegraphics[width=0.38\textwidth]{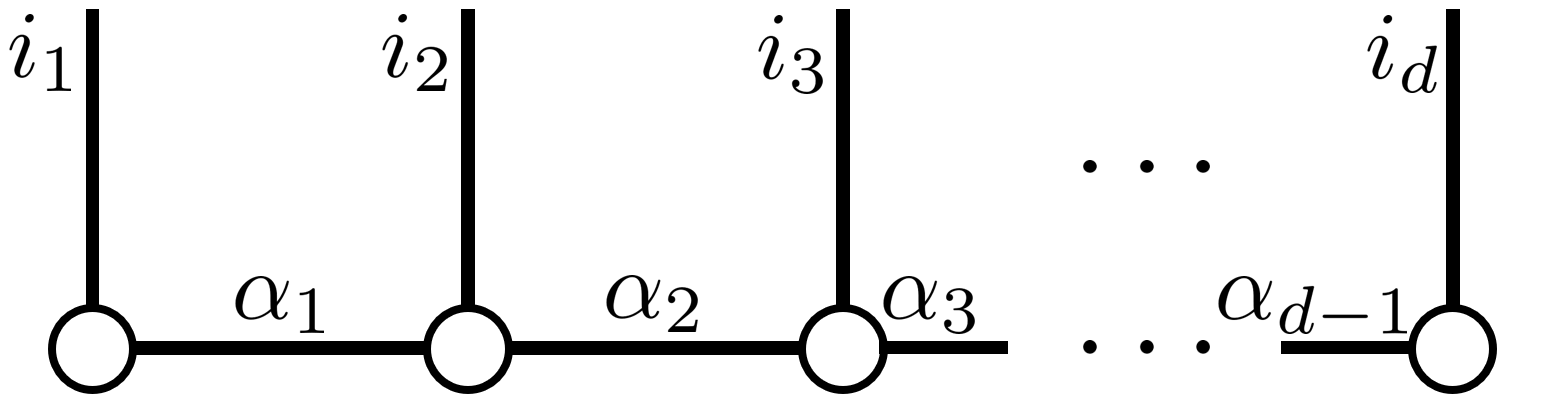}}} \quad \quad \quad \quad
    \subfloat[]{{\includegraphics[width=0.38\textwidth]{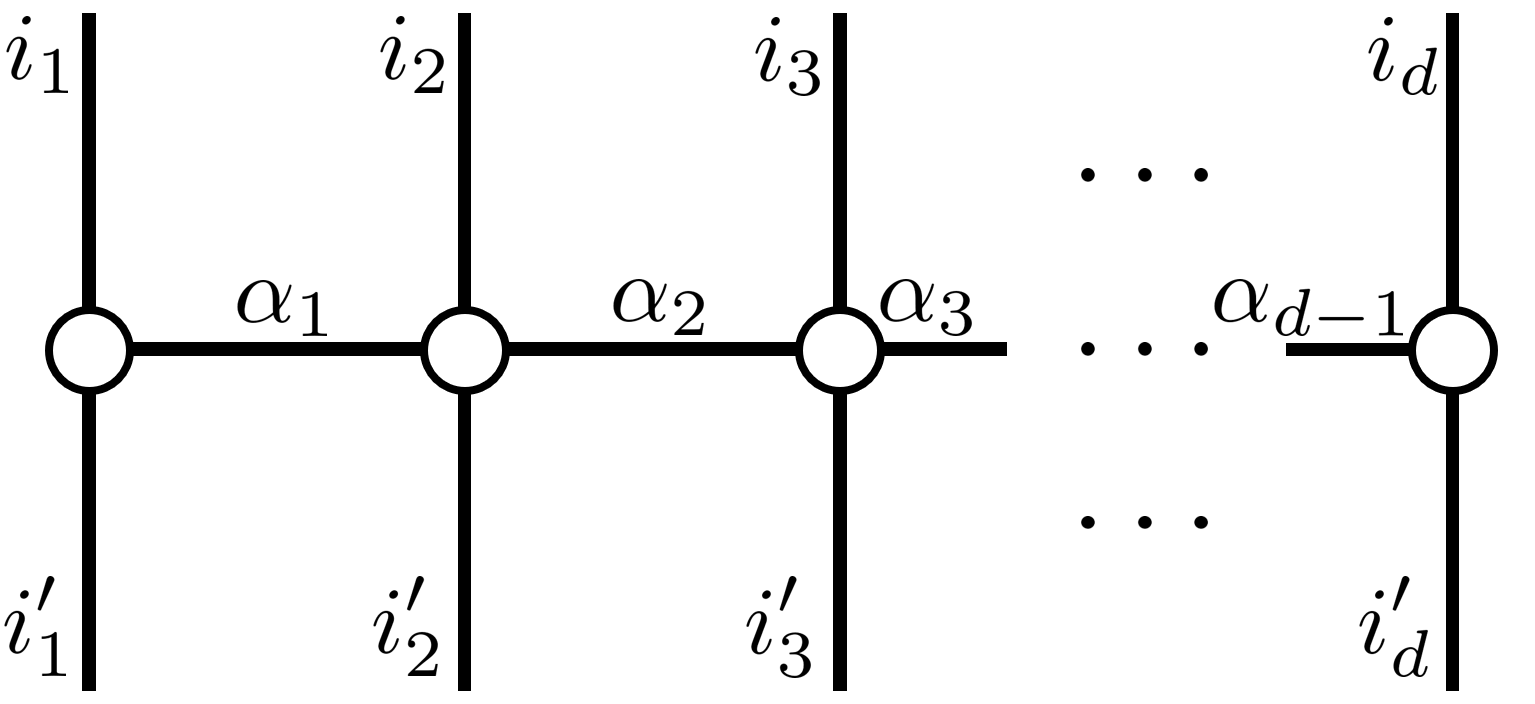}}}
    \caption{(a) $d$-dimensional TT/MPS. (b) $d$-dimensional MPO. The indices depicted match the expressions in \eqref{eq:TT} and \eqref{eq:MPO}. Note that we omit the legs for the trivial indices $\alpha_0$ and $\alpha_d$.}
    \label{fig:TT&FTT}
\end{figure}

In tensor diagrams, an MPS/TT is represented by a chain of $3$-tensors, as in \figref{fig:TT&FTT} (a). Note that the $0$-th and last tensor cores can be viewed as $2$-tensors since $r_0 = r_d = 1$, and as such the corresponding legs can be omitted from the diagram.

A matrix product operator (MPO) is a tensor network in which each constituent tensor has two external, uncontracted legs as well as two internal indices contracted with neighboring tensors in a chain-like fashion. Concretely, an MPO is a tensor
$\ccO \in\bbR^{(m_1\times m_2\times\dots\times m_d) \times (n_1\times n_2\times\dots\times n_d) }$
that can be written in the form
\begin{align}
    \ccO(i_1,\dots,i_d;i_1',\dots,i_d') \ = & \ \  \sum_{\alpha_0,\dots,\alpha_{d}} \ccG_1(\alpha_0,i_1,i_1',\alpha_1)\dots\ccG_d(\alpha_{d-1}, i_d,i_d',\alpha_d) \\ 
     = & \ \ \ccG_1(:,i_1,i_1',:)\dots\ccG_d(:, i_d,i_d',:),
    \label{eq:MPO}
\end{align}
where the sum over $\alpha_k$ has a range defined by a corresponding rank $r_k$, as in Definition \ref{def:mps}, and 
$r_0 = r_d = 1$ by convention. In this case we similarly say that the MPO has ranks $\vr=(r_0,\dots,r_{d})$. The $4$-tensor $\ccG_k \in \R^{r_{k-1} \times m_k \times n_k \times r_k}$ is called the $k$-th tensor core. A corresponding tensor network diagram is shown in \figref{fig:TT&FTT} (b).

For further background on tensor networks and diagrams, see \cite{tensor-diagram-intro}.

\section{Proposed method}
\label{sec:method}

In this work we obtain the committor function by solving the variational problem \eqref{eq:soft-varational} within a MPS/TT parametrization for the committor function $q$. We demonstrate that by approximating the equilibrium probability density $p$ in MPS/TT format, this optimization problem can be solved using basic tensor operations. In particular the minimization is accomplished using a standard alternating least squares approach. 

\subsection{Discretizing the variational problem}

We will represent the unknown committor function in a tensor product basis according to the product 
structure of the domain $\Omega=\Omega_1\times\Omega_2\times\dots\times\Omega_d$. Within this basis, we will approach the variational problem \eqref{eq:soft-varational} by Galerkin approximation.

\begin{figure}[!htb]
    \centering
    \subfloat{{\includegraphics[width=0.38\textwidth]{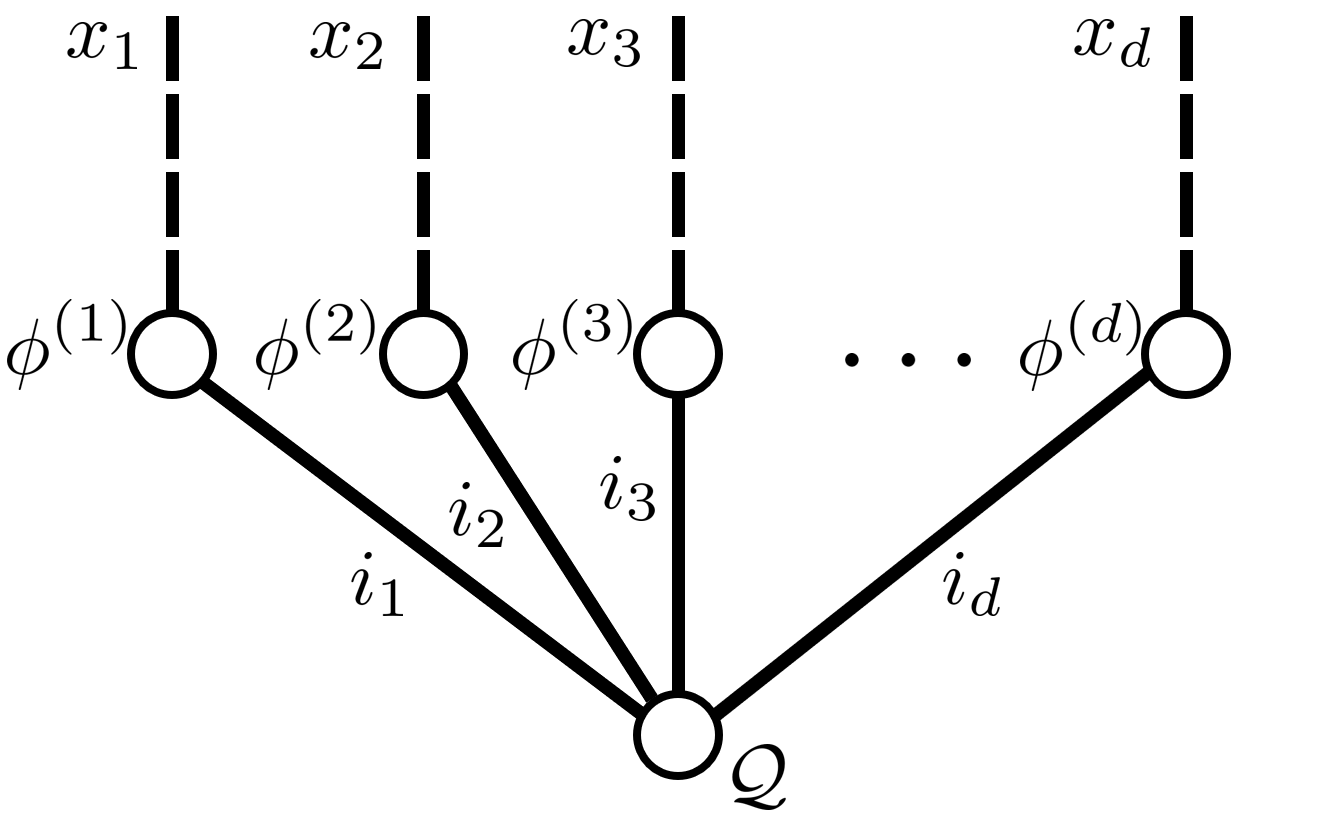}}}
    \caption{Tensor diagram for the decomposition  \eqref{eq:q-parametrization} of the committor function.}
    \label{fig:q}
\end{figure}

To begin, suppose that we have an orthogonal basis for each $L^2(\Omega_k)$, denoted by $\{\phi_j^{(k)} \}_{j=1}^{\infty}.$ In order to obtain a finite-dimensional problem we consider the subspace of $L^2(\Omega_k)$ spanned by only the first $L^{(k)}$ basis functions. Here the $L^{(k)}$, $k=1,\dots,d$, are a set of positive integers which are either fixed or determined adaptively. Then given the finite basis $\{\phi_j^{(k)}\}_{j=1}^{L^{(k)}}$, which spans a subspace of $L^2 (\Omega_k)$ for each $k=1,\dots,d$, we can consider 
an expansion of $q$ in the corresponding tensor product basis:
\begin{align}
    q(\vx) \ =& \ \  \sum_{i_1,\dots,i_d} \ccQ(i_1,\dots,i_d) \phi_{i_1}^{(1)}(x_1) \dots \phi_{i_d}^{(d)}(x_d), \nonumber\\
    :=&\ \  \sum_{\bvec{i}}\ccQ(\bvec{i})\phi_{i_1}^{(1)}(x_1) \dots \phi_{i_d}^{(d)}(x_d),
    \label{eq:q-parametrization}
\end{align}
where, for notational convenience, we have defined $\bvec{i}:=(i_1,i_2,\dots,i_d)$ and $\vx:=(x_1,x_2,\dots,x_d)$. Additionally, we set $\phi^{(k)}:=( \phi_{j}^{(k)} )_{j=1}^{L^{(k)}}$. Each $\phi^{(k)}$ can be viewed as a 2-tensor via $\phi^{(k)} (j, x) = \phi^{(k)}_{j} (x)$, where the first index is the basis function index and the second (continuous) index is a spatial coordinate. Then the decomposition \eqref{eq:q-parametrization} for $q$ can be depicted graphically as in  \figref{fig:q}.

To determine the committor function $q$ within the truncated tensor product basis, we want to determine the coefficient tensor $\ccQ$ that is optimal in the sense of \eqref{eq:soft-varational}. 
By inserting the parameterization \eqref{eq:q-parametrization} into the variational problem \eqref{eq:soft-varational}, we can rewrite the optimization problem as follows,
\begin{align}
    \argmin_{\ccQ} & \quad  \underbrace{\sum_{k=1}^d \sum_{\vi,\vj} H^k(\vi;\vj) \ccQ(\vi)\ccQ(\vj)}_{\int_{\Omega} \abs{\nabla q(\vx)}^2 p(\vx)\, d\vx} + \underbrace{\rho \sum_{\vi,\vj} H^A(\vi;\vj) \ccQ(\vi)\ccQ(\vj)}_{\rho \int_{\Omega} q(\vx)^2 p_{A}(\vx)\, d\vx} \nonumber\\
    & \quad \quad \quad + \  \underbrace{\rho \sum_{\vi,\vj} H^B(\vi;\vj) \ccQ(\vi)\ccQ(\vj) - 2\rho \sum_{\vi} \ccQ(\vi) h^B(\vi) + \rho}_{\rho \int_{\Omega} (q(\vx) - 1)^2 p_{B}(\vx) \,d\vx},
    \label{eq:variational-tensor}
\end{align}
where 
\begin{align}
    & H^k(\vi;\vj) = \int_{\Omega} \pfrac{}{x_k} \left[\phi_{i_1}^{(1)}(x_1) \dots \phi_{i_d}^{(d)}(x_d)\right] \pfrac{}{x_k} \left[\phi_{j_1}^{(1)}(x_1) \dots \phi_{j_d}^{(d)}(x_d)\right] p(\vx) \, d\vx 
    \label{eq:H^k}\\
    & H^A(\vi;\vj) = \int_{\Omega} \phi_{i_1}^{(1)}(x_1) \dots \phi_{i_d}^{(d)}(x_d) \phi_{j_1}^{(1)}(x_1) \dots \phi_{j_d}^{(d)}(x_d) p_{A}(\vx) \,  d\vx \label{eq:H^A}\\
    & H^B(\vi;\vj) = \int_{\Omega} \phi_{i_1}^{(1)}(x_1) \dots \phi_{i_d}^{(d)}(x_d) \phi_{j_1}^{(1)}(x_1) \dots \phi_{j_d}^{(d)}(x_d) p_{B}(\vx) \, d\vx \label{eq:H^B}\\
    & h^B(\vi) = \int_{\Omega} \phi_{i_1}^{(1)}(x_1) \phi_{i_2}^{(2)}(x_2)\dots \phi_{i_d}^{(d)}(x_d) p_{B}(\vx) \,  d\vx. \label{eq:h^B}
\end{align}

Here $(\vi;\vj)=(i_1,\dots,i_d;j_1,\dots,j_d)$ is a concatenation of multi-indices. We can simply ignore the last constant $\rho$ since it does not affect the minimizer. Computing the tensors $\{H^k\}_{k=1}^d$, $H^A$, $H^B$, and $h^B$ is \emph{prima facie} intractable as it requires us to perform integration over the $d$-dimensional domain $\Omega$, in addition to storing tensors of exponential size in $d$. Moreover, the number of unknown tensor entries of $\ccQ$ is also exponential in $d$. Traditional approaches are therefore prohibitively expensive for $d$ of even moderate size.

In the next two sections we show how to use MPS/TT approximations to obtain $\{H^k\}_{k=1}^d,H^A,H^B,h^B,$ allowing us to solve the optimization problem \eqref{eq:variational-tensor} with computational and storage complexities of $\ccO(d)$.

\subsection{Constructing $H^k$}
\label{sec:construc-H^k}

\begin{figure}[!htb]
    \centering
    \subfloat{{\includegraphics[width=0.47\textwidth]{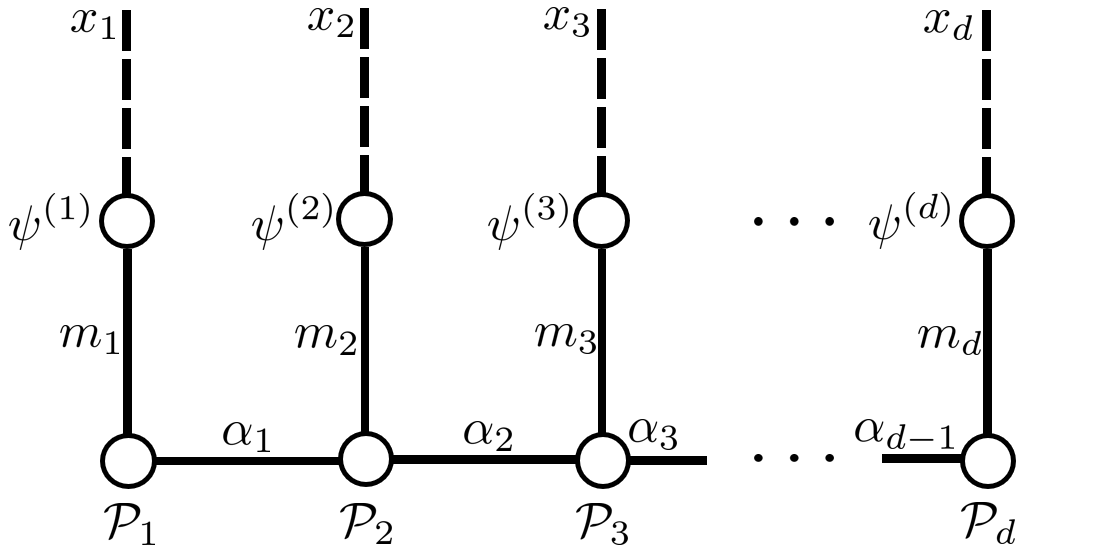}}}
    \caption{Tensor diagram for the approximate equilibrium density $p$ \eqref{eq:p-FTT}.}
    \label{fig:p}
\end{figure}

In this subsection we detail the construction of $H^k$, which corresponds to the variational energy term $\int_{\Omega} \abs{\nabla q(\vx)}^2 p(\vx) d\vx$ of \eqref{eq:soft-varational}. As mentioned above, in order to obtain each $H^k$ in \eqref{eq:H^k}, one needs to evaluate a $d$-dimensional integral and store the resulting high-dimensional tensor. To circumvent the exponential complexity in $d$, we assume that the equilibrium density $p$ can be approximated as an MPS/TT as follows:
\begin{align}
    p(\vx) = \sum_{\substack{m_1,\dots,m_d\\\alpha_0,\dots,\alpha_{d}}} &\ccP_1(\alpha_0,m_1,\alpha_1)\dots\ccP_d(\alpha_{d-1}, m_d, \alpha_d) \psi_{m_1}^{(1)}(x_1) \cdots \psi_{m_d}^{(d)}(x_d),
    \label{eq:p-FTT}
\end{align}
where $\psi^{(k)}:=( \psi_{j}^{(k)} )_{j=1}^{K^{(k)}}$ is a vector of univariate basis functions $\Omega_i\rightarrow\bbR$. \figref{fig:p} illustrates the structure of the equilibrium density $p$ that we assume in this paper.

The construction of the MPS/TT format for a given equilibrium density $p$ will be described in Section \ref{sec:numerical-tests} in the contexts of specific example problems.

\begin{figure}[!htb]
    \centering
    \subfloat[$H^k$]{{\includegraphics[height=0.4\textwidth]{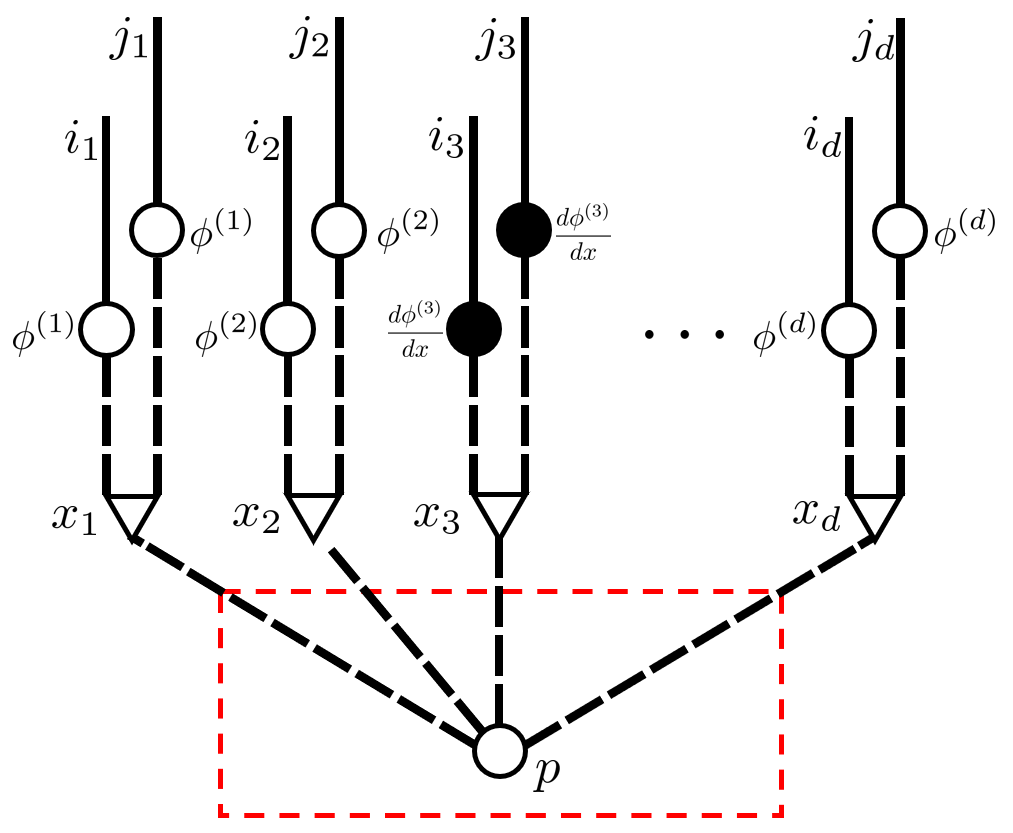}}} \quad \quad \quad
    \subfloat[Approximation of $H^k$]{{\includegraphics[height=0.4\textwidth]{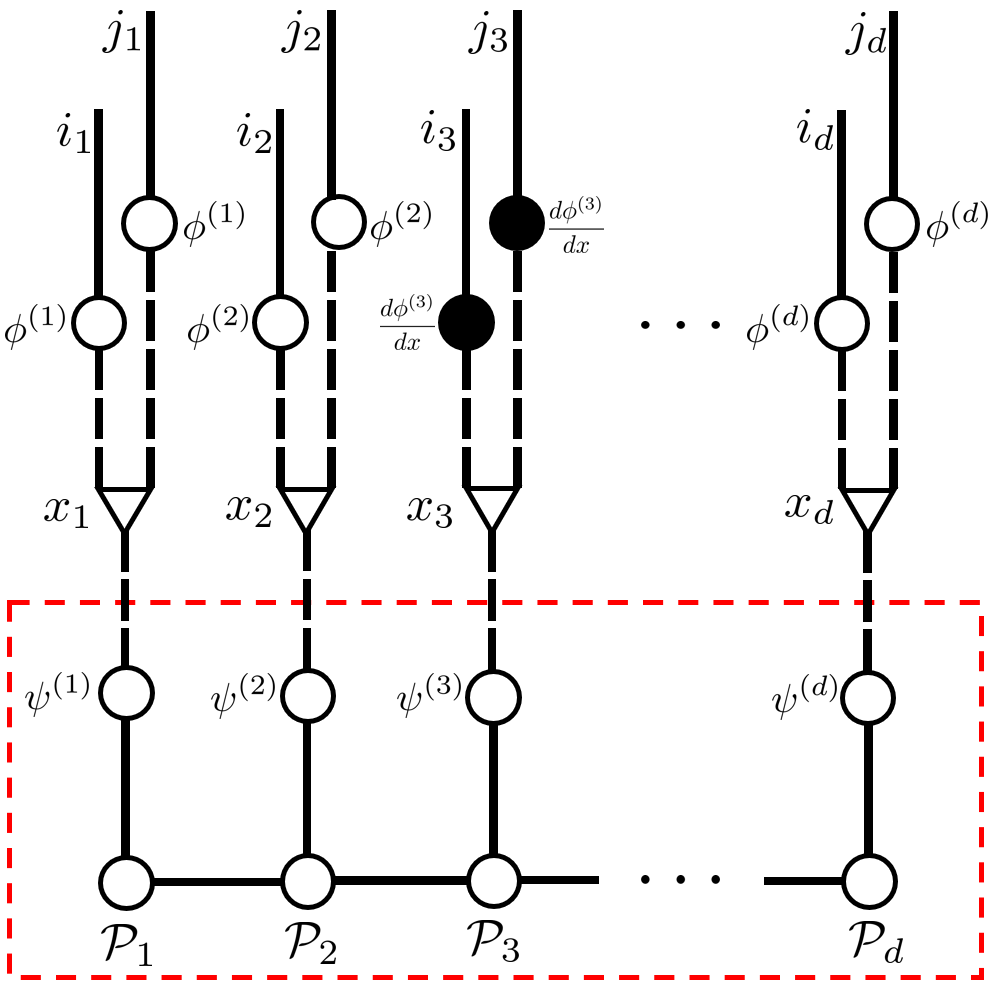}}}
    \caption{(a) Tensor diagram for $H^k$, $k=3$, as defined in \eqref{eq:H^k}. (b) Approximation of $H^k$ obtained by replacing $p$ with its MPS/TT approximation \eqref{eq:p-FTT}. We use red dashed boxes to indicate the region of replacement. The original basis functions $\phi^{(l)}$, $l \neq k$, are represented using hollow nodes, and the derivative $d\phi^{(k)}/dx$ is distinguished using a filled node.}
    \label{fig:H^k-original-vs-simplified}
\end{figure}

\begin{figure}[!htb]
    \centering
    \subfloat[$H_l$]{{\includegraphics[width=0.35\textwidth]{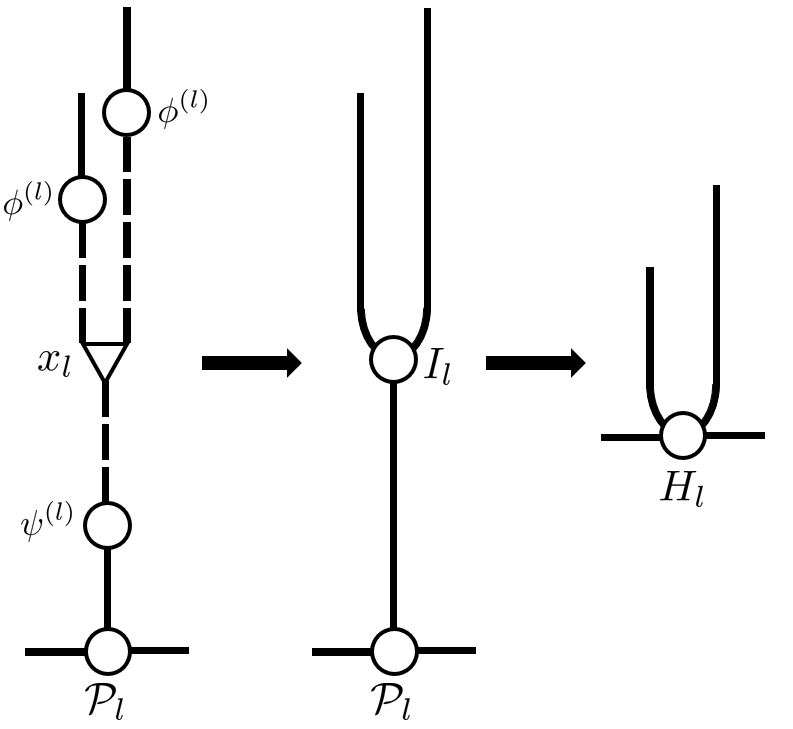}}}\quad\quad \quad \quad 
    \subfloat[$\tilde H_l$]{{\includegraphics[width=0.35\textwidth]{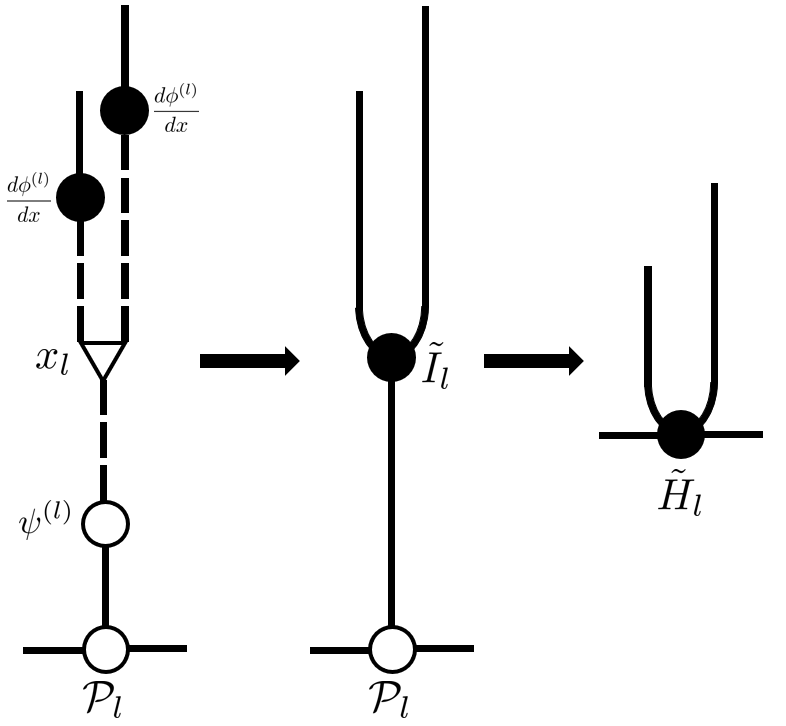}}}\\ \vspace{2 mm}
    \subfloat[]{{\includegraphics[width=0.75\textwidth]{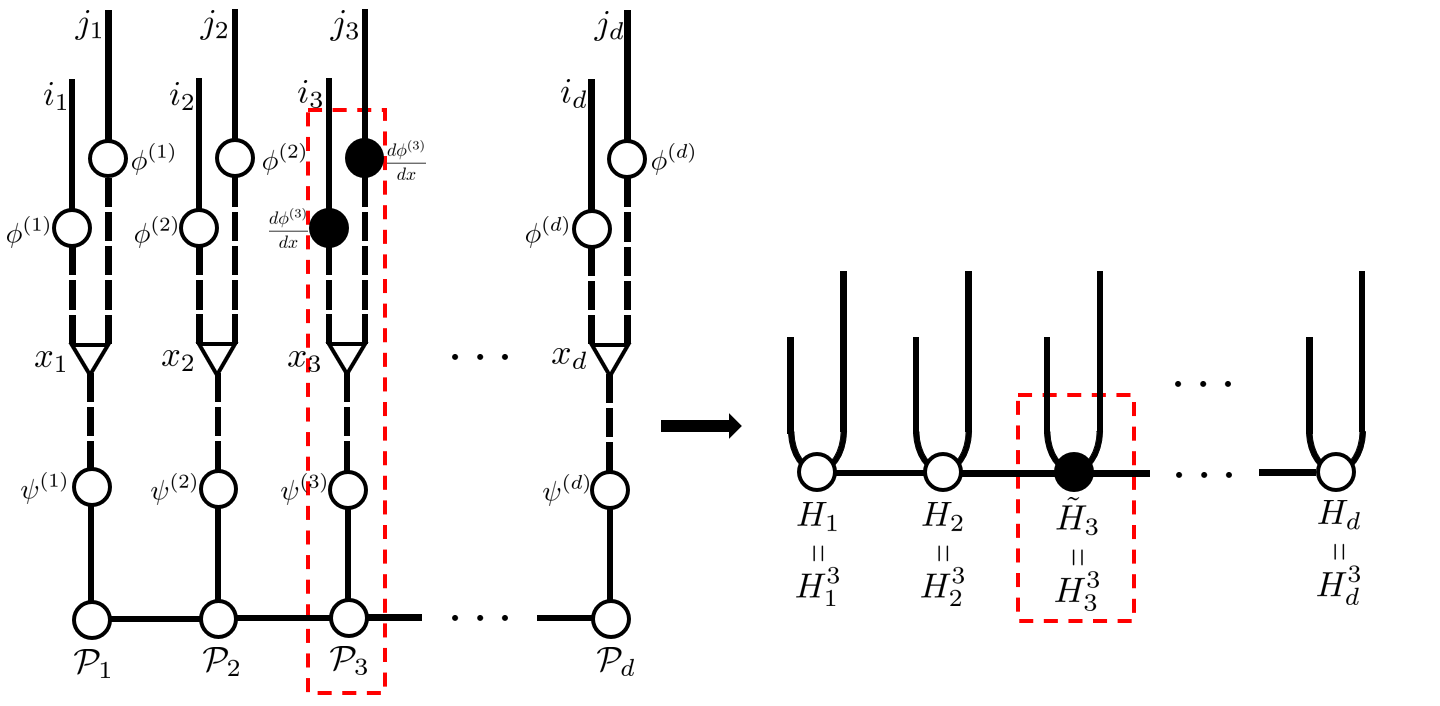}}}
    \caption{Diagrammatic illustration of the construction of an MPO format for $H^k$. (a) Precompute tensors $H_l$. (b) Precompute tensors $\tilde H_l$. (c) Assemble $H^k$, $k=3$, by substituting tensors $H_l$ and $\tilde H_l$. We use red dashed boxes to indicate the region of replacement. A similar construction is repeated for other $k$.}
    \label{fig:var-int}
\end{figure}

Such an approximation of $p$ amounts to changing the tensor representation of $H^k$ depicted graphically in \figref{fig:H^k-original-vs-simplified} (a) to the representation in \figref{fig:H^k-original-vs-simplified} (b).
Note that these calculations involve the derivatives of our univariate basis functions $\phi^{(k)}$. In our figures, we use a hollow node to represent the vector of basis functions $\phi^{(k)}$ and a filled node to represent its vector of derivatives, i.e., $d\phi^{(k)}/dx$.
We observe that the $\{H^k\}_{k=1}^{d}$ are naturally viewed as MPOs, and moreover the construction of these MPOs can be performed using basic tensor algebra in $\ccO(d)$ complexity. More precisely:
~\\
\begin{enumerate}[{(1)}]
    \item To construct an MPO for $H^k$ following \figref{fig:H^k-original-vs-simplified} (b), note that we need to perform two types of tensor contraction: one involving the original basis functions $\phi^{(l)}$, $l\neq k$, and the other involving the derivatives $d\phi^{(k)}/dx$. Therefore we precompute these two contractions, which can be recycled to form MPOs for the $H^k$, $k=1,\ldots ,d$. First, we form $I_l$, $l=1,\ldots,d$ by contracting three tensors $\psi^{(l)}$, $\phi^{(l)}$, $\phi^{(l)}$ and form $\tilde I_l$, $l=1,\ldots,d$ by contracting three tensors $\psi^{(l)}$, $d\phi^{(l)}/dx$, $d\phi^{(l)}/dx$. The tensors $I_l$ and $\tilde{I}_l$ are defined graphically in \figref{fig:var-int} (a) and (b). These contractions can be performed by univariate numerical integration. Next we contract each $I_l$ with the corresponding tensor core $\ccP_l$ to obtain $H_l$ for $l=1,\dots,d$. Similarly we contract $\tilde I_l$ with $\ccP_l$ to obtain $\tilde H_l$. These constructions are illustrated in \figref{fig:var-int} (a) and (b), respectively.
    
    \item Next we assemble $H^k$ by substituting $H_l,\tilde H_l$, $l=1,\ldots,d$ into the red boxes in \figref{fig:var-int} (c) as needed. This yields an MPO as shown in \figref{fig:var-int} (c) on the right. We denote $l$-th tensor core of $H^k$ by $H_l^k$.
    
    \item Finally we repeat step (2) for all $k=1,\dots,d$.
\end{enumerate}
~\\

The algorithm outputs core tensors $H_1^k,\dots ,H_d^k$ for the MPO $H^k$, $k=1,\dots,d$. The computational complexity of forming each pair $I_l,\tilde I_l$ is $\ccO(K^{(l)} L^{(l)2})$ since a univariate numerical integration is performed for each entry. Hence the cost of computing all of the $I_l,\tilde I_l$ is $\ccO(d K L^2)$, where we define 
$K := \max_l K^{(l)}$, $L := \max_l L^{(l)}$.

If we assume that the MPS/TT format for $p$ has ranks $\vr = (r_0,r_1,r_2,\dots,r_d)$, and set $r = \max_l r_l,$ then the contraction steps with the $\ccP_l$ altogether cost $\ccO(d r^2 K L^2)$. Therefore to construct the tensor cores $H_1^k,\dots,H_d^k$, the total computational cost is $\ccO(d r^2 K L^2)$. The memory complexity, including that of storing the intermediate tensors $I_l$, $\tilde I_l$, $H_l$, $\tilde H_l$, is $\ccO(d KL^2 + d r^2 L^2)$.

\subsection{Constructing $H^A$ and $H^B$}
\label{sec:construct-H^A-H^B}

The tensors $H^A$ and $H^B$ derive from the two soft boundary penalty terms $\rho \int_{\Omega} q(\vx)^2 p_{A}(\vx) \,d\vx$ and $\rho \int_{\Omega} (q(\vx) - 1)^2 p_{B}(\vx) \, d\vx$ in \eqref{eq:soft-varational}. The construction of the MPO format for $H^A$ and $H^B$ is similar to that of $H^k$ detailed above in \secref{sec:construc-H^k}. However, 
we now further need to represent the soft boundary measures $p_A$ and $p_B$ as MPS/TT. Recall our motivation that $p_A$ and $p_B$ weakly approximate surface measures on the boundaries $\partial A$ and $\partial B$, though for \emph{any} choice of $p_A$ and $p_B$ we may still interpret the soft committor function probabilistically following Appendix \ref{app:soft}. The construction of approximate surface measures varies depending on the specific geometry of $A$ and $B$. In many applications, $A$ and $B$ are balls or half-spaces, so $\partial A$ and $\partial B$ are spheres or hyperplanes. We discuss these two cases in detail presently.

For the case of a sphere, the Gaussian annulus theorem \cite[Theorem 2.8]{foundation-of-ds} indicates that most of the mass of a high-dimensional Gaussian distribution concentrates on a shell. If we assume that region $A$ is a $d$-dimensional ball with center $\vx_A$ and radius $R_A$, we can approximate the uniform measure on $\partial A$ by a Gaussian density,
\begin{align}
    p_{A}(\vx) = \frac{1}{(2\pi)^{d/2} \sigma^d} \exp{\left(-\frac{\|\vx - \vx_A\|^2}{2\sigma^2}\right)}, \ \text{where}\ \ \sigma=\frac{R_A}{\sqrt{d}}.
    \label{eq:high-d-gaussian-for-sphere}
\end{align}
More precisely, under this probability measure, we have that
\begin{equation}
\text{Prob}\left(\left \vert \|\vx\|_2 - \sqrt d \sigma \right\vert \geq t\sigma\right)\leq \exp\left(-\frac{t^2}{\kappa} \right),\ \text{for all}\ t>0,
\end{equation}
where $\kappa>0$ is a constant \cite{vershynin2018high}. This bound indicates that the mass of $p_A$ concentrates on a shell with radius $R_A$ and thickness $O(1/\sqrt d)$.  It is  straightforward to convert \eqref{eq:high-d-gaussian-for-sphere} into MPS/TT format since it is in fact a pure tensor product of univariate functions of each scalar variable $x_k$. As such our resulting MPS/TT should have ranks all equal to 1.

For the case of a hyperplane, suppose in particular that $\partial A = \{\vx \in \R^d \,:\,x_i=c\}$. In this case consider
\begin{align}
    p_{A}(\vx) = \frac{1}{(2\pi)^{1/2} \sigma} \exp{\left(-\frac{(x_i - c)^2}{2\sigma^2}\right)},
    \label{eq:gaussian-for-hyperplane}
\end{align}
where $\sigma$ controls the sharpness of the approximation. (Note that $p_A$ is not a normalized probability density because $\partial A$ is not compact, but in fact the interpretation of Appendix \ref{app:soft} applies to an arbitrary nonnegative function $p_A$, not necessarily integrable.)
This choice of $p_A$ is again a pure tensor product of univariate functions.

\begin{figure}[!htb]
    \centering
    \subfloat{{\includegraphics[width=0.4\textwidth]{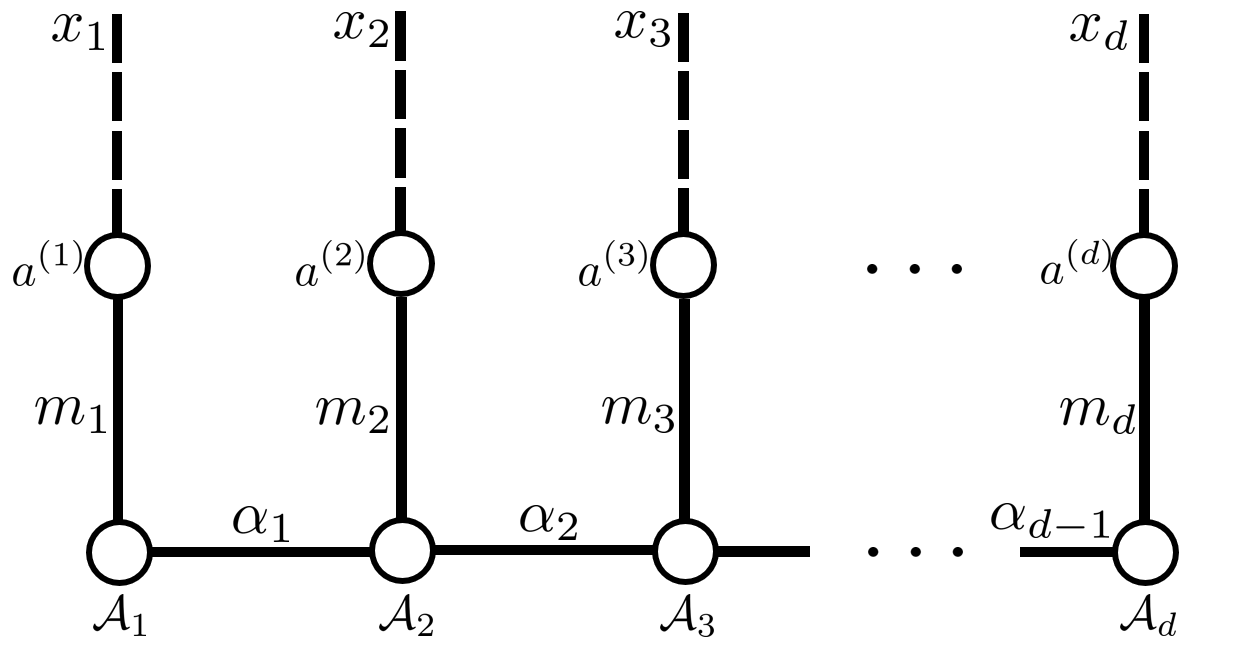}}} \quad \quad \quad \quad
    \subfloat{{\includegraphics[width=0.4\textwidth]{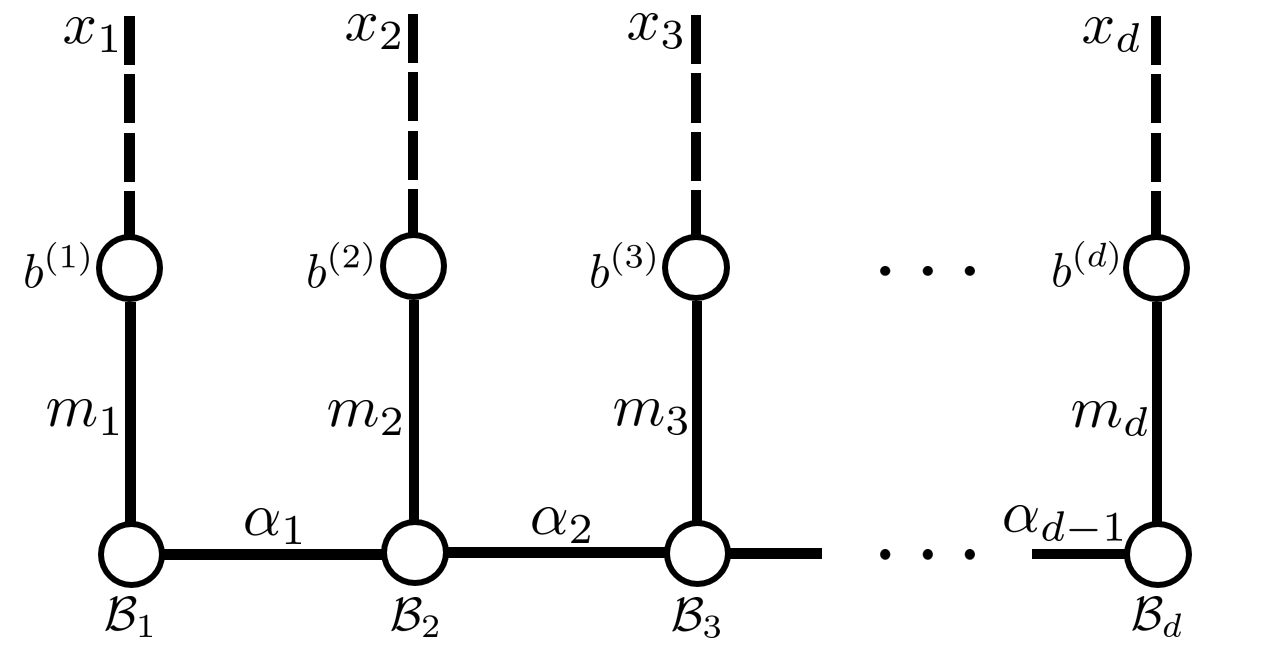}}}
    \caption{Tensor diagrams for (a) $p_{A}$ as in \eqref{eq:p_A-FTT} and (b) $p_{B}$ as in \eqref{eq:p_B-FTT}.}
    \label{fig:p_B}
\end{figure}

Now suppose that we have $p_{A}$ and $p_{B}$ in MPS/TT form. Specifically, assume that we have tensor cores $\{\ccA_i\}_{i=1}^d$ with associated ranks $\bm{s} = (s_0,s_1,\dots,s_d)$, together with a vector of basis functions $a^{(k)} := ( a_j^{(k)} )_{j=1}^{K_A^{(k)}}$ for each $k=1,\ldots, d$. Also assume that we have tensor cores $\{\ccB_i\}_{i=1}^d$ with ranks $\bm{t} = (t_0,t_1,\dots,t_d)$ and a vector of basis functions $b^{(k)} := ( b_j^{(k)} )_{j=1}^{K_B^{(k)}}$. Then we assume that we can write $p_A$ and $p_B$ as 

\begin{align}
    & p_{A}(\vx) = \sum_{\substack{m_1,\dots,m_d\\\alpha_0,\dots,\alpha_{d}}} &\ccA_1(\alpha_0,m_1,\alpha_1)\dots\ccA_d(\alpha_{d-1}, m_d, \alpha_d) a_{m_1}^{(1)}(x_1) \dots a_{m_d}^{(d)}(x_d).
    \label{eq:p_A-FTT}\\
    & p_{B}(\vx) = \sum_{\substack{m_1,\dots,m_d\\\alpha_0,\dots,\alpha_{d}}} &\ccB_1(\alpha_0,m_1,\alpha_1)\dots\ccB_d(\alpha_{d-1}, m_d, \alpha_d) b_{m_1}^{(1)}(x_1) \dots b_{m_d}^{(d)}(x_d).
    \label{eq:p_B-FTT}
\end{align}
In \figref{fig:p_B} we illustrate these formats graphically.

Now in light of the resemblance among \eqref{eq:H^k}, \eqref{eq:H^A},and \eqref{eq:H^B}, we can use the same procedure described in \secref{sec:construc-H^k} to approximate $H^A$ and $H^B$ as MPOs. To wit, we simply replace $p$ in \figref{fig:H^k-original-vs-simplified} with the MPS/TT approximations of $p_{\partial A}$ or $p_{\partial B}$ and replace all derivatives $d\phi^{(l)}/dx$ by $\phi^{(l)}$ since there are no derivatives in \eqref{eq:H^A} and \eqref{eq:H^B}. Ultimately we obtain MPO formats for $H_A$ and $H_B$ with ranks $\bm{s}$ and $\bm{t}$ and cores $H_k^{A}$ annd $H_k^{B}$, $k=1,\ldots, d$ respectively. The computational complexities of constructing $H^A$ and $H^B$ are $\ccO(d s^2 K_A L^2)$ and $\ccO(d t^2 K_B L^2)$, respectively, where we define $s := \max_l s_l$, $t := \max_l t_l$, $K_A := \max_l K_A^{(l)}$, and $K_B := \max_l K_B^{(l)}$. The memory complexities are $\ccO(d K_A L^2 + d s^2 L^2)$ and $\ccO(d K_B L^2 + d t^2 L^2)$, respectively.

\subsection{Constructing $h^B$}
\label{sec:construct-h^B}

\begin{figure}[!htb]
    \centering
    \subfloat[$h^B$]{{\includegraphics[height=0.28\textwidth]{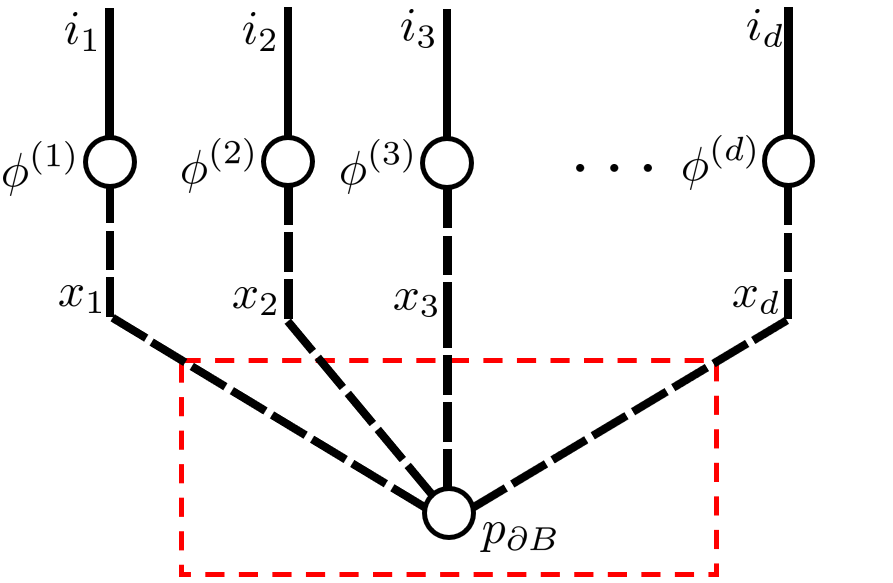}}} \quad \quad \quad
    \subfloat[Approximation of $h^B$]{{\includegraphics[height=0.28\textwidth]{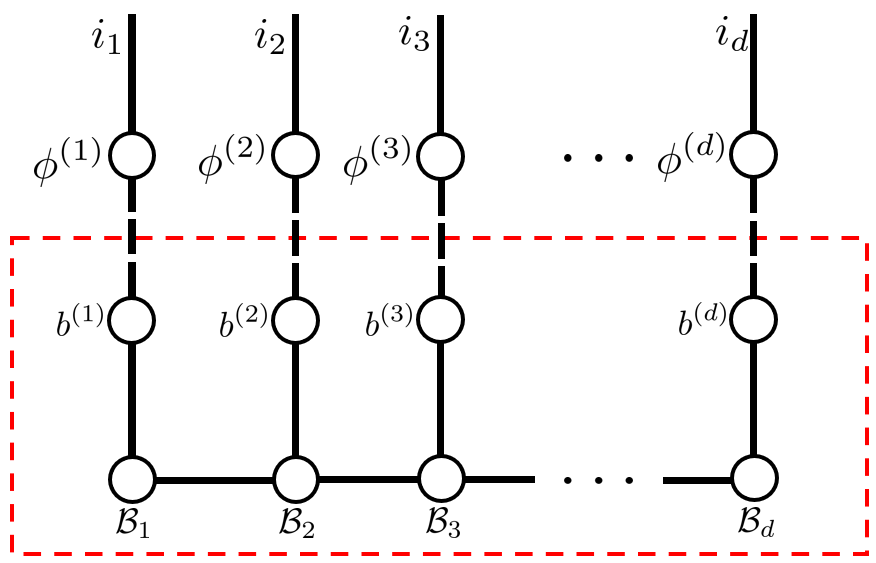}}}
    \caption{(a) Tensor diagram for $h^B$ as in \eqref{eq:h^B}. (b) Approximation of $h^B$ obtained by replacing the soft boundary measure $p_{B}$ with its MPS/TT approximation \eqref{eq:p_B-FTT}. We use red dashed boxes to indicate the region of replacement.}
    \label{fig:h^B-original-vs-simplified}
\end{figure}

\begin{figure}[!htb]
    \centering
    \subfloat[]{{\includegraphics[width=0.8\textwidth]{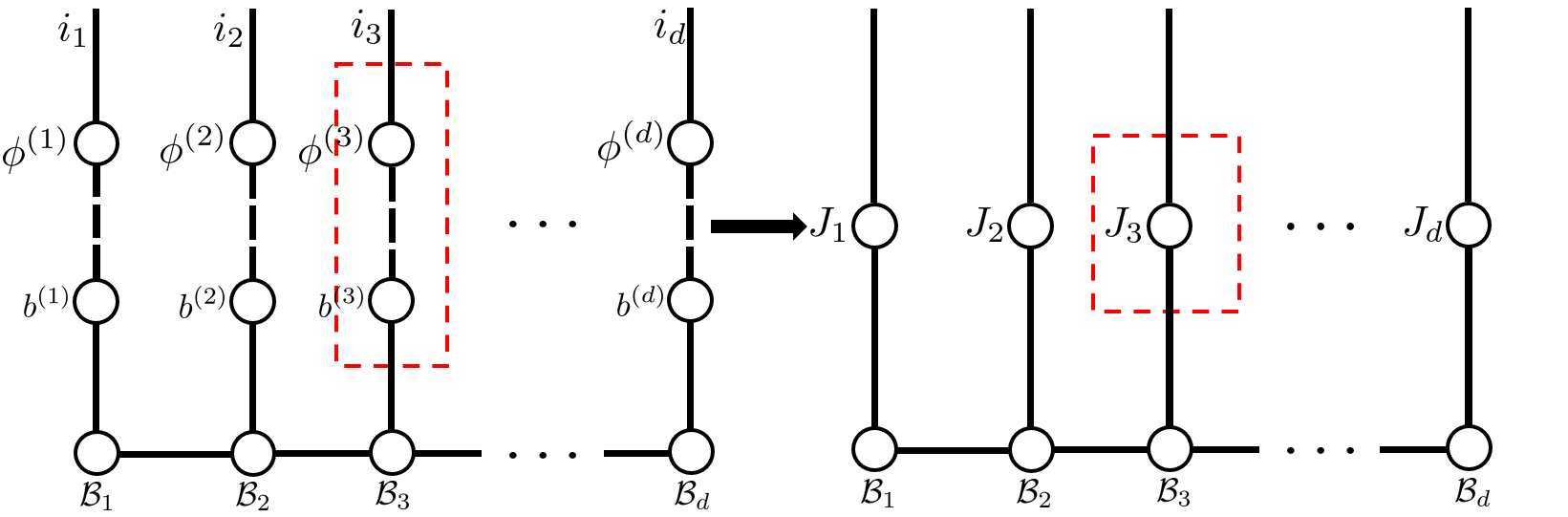}}}\\ \vspace{2 mm}
    \subfloat[]{{\includegraphics[width=0.8\textwidth]{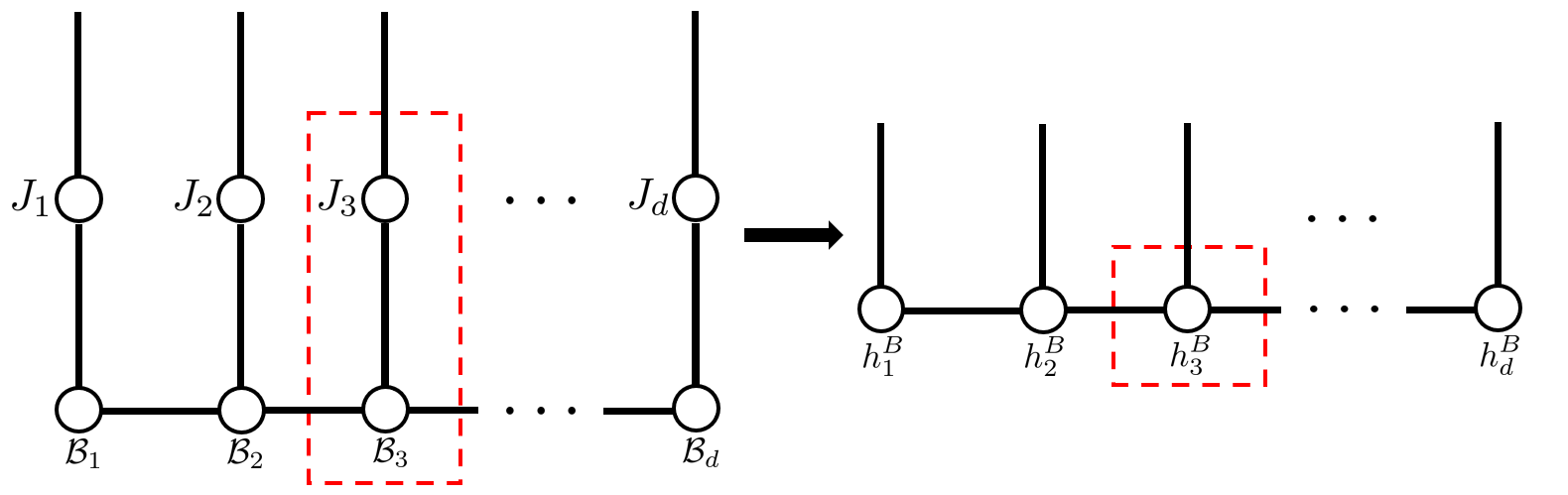}}}
    \caption{Illustration of the construction of $h^B$ in MPS/TT format. (a) Contract tensors $\phi^{(l)}$ and $b^{(l)}$ to get $J_l$. (b) Contract tensors $J_l$ and $\ccB_l$ to get $h_l^B$. We use red dashed boxes emphasize the contractions.}
    \label{fig:bnd-int}
\end{figure}

In this subsection we focus on constructing $h^B$, which comes from the cross term in the second penalty term $\rho \int_{\Omega} (q(\vx) - 1)^2 p_{B}(\vx) \, d\vx$ within \eqref{eq:soft-varational}. The ideas are again very similar to \secref{sec:construc-H^k}. The tensor diagram for $h^B$ is shown in \figref{fig:h^B-original-vs-simplified} (a). By plugging in the MPS/TT approximation of the soft boundary measure $p_{B}$ \eqref{eq:p_B-FTT}, we obtain the approximation of $h^B$ illustrated in \figref{fig:h^B-original-vs-simplified} (b). One can further bring $h^B$ to a standard MPS/TT form, using the contractions shown in \figref{fig:bnd-int}. In detail, the procedure is as follows:
~\\
\begin{enumerate}[{(1)}]
    \item In \figref{fig:bnd-int} (a), we contract the two connected basis function nodes in the red box. This results in tensors $J_k,\ k=1,\ldots,d$, seen in \figref{fig:bnd-int} (a) on the right. This contraction requires univariate numerical integrations. 
    
    \item Next we merge the computed tensor $J_k$ and the tensor core $\ccB_k$ for $k=1,\ldots,d$. The resulting 3-tensors, which are the tensor cores for $h^B$, are denoted $h_k^B$. This step is shown in \figref{fig:bnd-int} (b).
\end{enumerate}
~\\

This procedure yields $h^B$ in MPS/TT format with tensor cores $h_1^B,\dots ,h_d^B.$ The computational and memory complexities of constructing $h^B$ are $\ccO(d t^2 K_B L)$ and $\ccO( d K_B L + d t^2 L )$, respectively.

\subsection{Optimization}
\label{sec:optimization}

We have discussed how the MPS/TT format can be used to compress the tensors $\{H^k\}_{k=1}^d$, $H^A$, $H^B$, and $h^B.$ In order to obtain a tractable algorithm for computing the committor function, it is natural to represent the unknown tensor $\ccQ$ in a compatible format. Indeed, without imposing some additional structure on the parameterization \eqref{eq:q-parametrization}, the unknown tensor core $\ccQ$ is still of size exponential in $d$. Thus we approximate $\ccQ$ as in MPS/TT format as 
\begin{align}
    \ccQ(\vi) := \sum_{\alpha_0,\dots,\alpha_{d}} \ccQ_1(\alpha_0,i_1,\alpha_1)\ccQ_2(\alpha_1, i_2, \alpha_2)\dots\ccQ_d(\alpha_{d-1}, i_d, \alpha_d).
    \label{eq:Q-TT}
\end{align}

\begin{figure}[!htb]
    \centering
    \subfloat{{\includegraphics[width=0.45\textwidth]{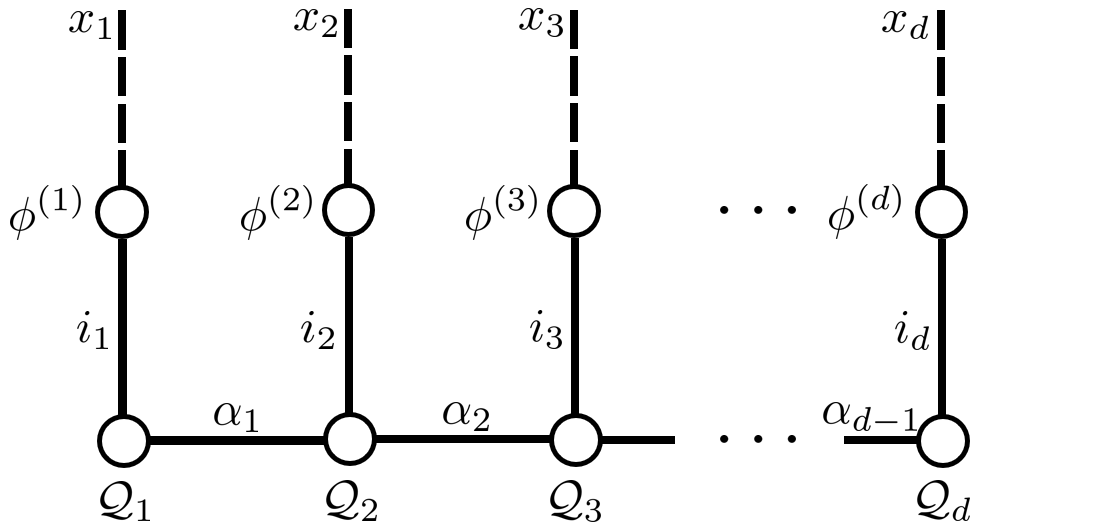}}}
    \caption{Tensor diagram for the parametrization of $q$ following \eqref{eq:q-parametrization} and \eqref{eq:Q-TT}.}
    \label{fig:q-FTT}
\end{figure}

\begin{figure}[!htb]
    \centering
    \subfloat[Original variational problem]{{\includegraphics[width=1\textwidth]{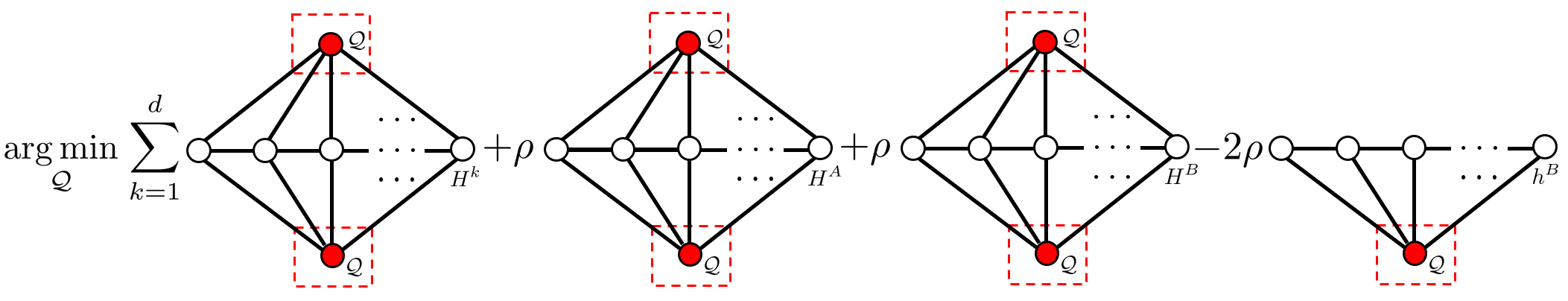}}}\\ \vspace{2 mm}
    \subfloat[Approximated variational problem]{{\includegraphics[width=1\textwidth]{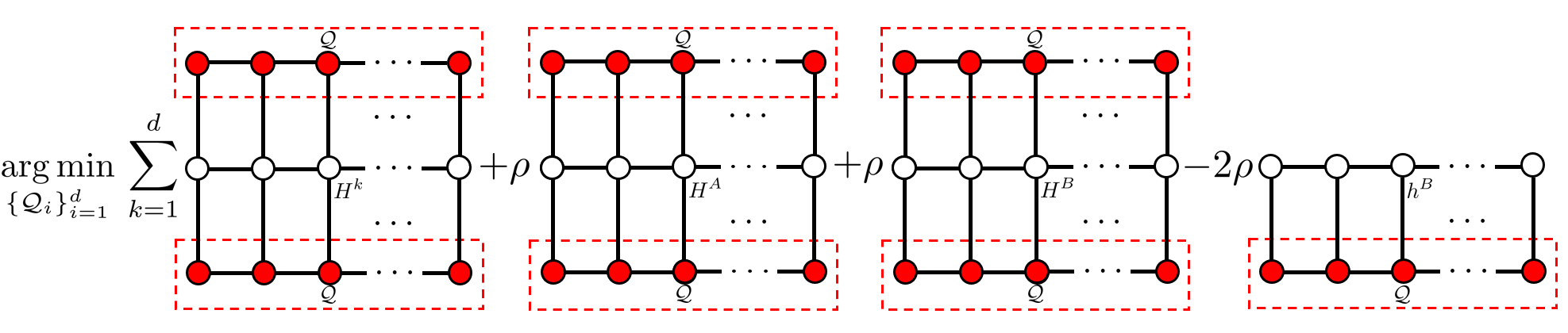}}}
    \caption{(a) Tensor diagram representation of the variational problem corresponding to \eqref{eq:variational-tensor}. (b) Approximate variational problem obtained by replacing $\ccQ$ with its MPS/TT approximation \eqref{eq:Q-TT}. We use red dashed boxes to indicate the region of replacement. Specifically the red tensor cores are unknown variables in the optimization.}
    \label{fig:var-problem}
\end{figure}

The tensor diagram for $q$ is shown in \figref{fig:q-FTT}. Empirically we observe that this format is able to capture the structure of $q$ accurately, i.e., without growth of the ranks of the tensor cores. The MPS/TT format \eqref{eq:Q-TT} for $q$ greatly simplifies the solution of the variational problem \eqref{eq:variational-tensor}. In \figref{fig:var-problem} we compare the tensor diagram depictions of the original variational problem and the new simplified problem by replacing $\ccQ$ with its MPS/TT approximation. We note that all terms in the simplified form (\figref{fig:var-problem} (b)) can be computed with standard MPO-MPS or MPS-MPS contractions in $\ccO(d)$ time.

In \figref{fig:var-problem}, the unknown tensor cores of MPS $\ccQ$ are marked in red. A standard approach for optimization problems of the form \figref{fig:var-problem} (b) is alternating least squares (ALS). In each ALS iteration, we loop over the dimensions $k=1,\dots,d$. For each $k$, we treat all coefficient tensor cores but $\ccQ_k$  as constant. Thie yields an unconstrained least squares problem for $\ccQ_k$. Naively the computational complexity is $\ccO(d^2)$ since the bottleneck is the summation of $d$ terms in \figref{fig:var-problem} (b) and each term requires at least $\ccO(d)$ tensor contractions. However by using the same trick as in the construction of $H^k$ and carefully reusing the computed nodes, one can bring the computational complexity down to $\ccO(d)$.

\section{Numerical experiments}
\label{sec:numerical-tests}

In this section, we present numerical results that demonstrate the accuracy and efficiency of the proposed method. 

\subsection{Double-well potential}

In the first numerical experiment, we consider the following potential
\begin{align}
    V(\vx) = (x_1^2 - 1)^2 + 0.3\sum_{i=2}^d x_i^2, 
    \label{eq:DW-potential}
\end{align}
and we let $A$, $B$ be the half-spaces
\begin{align}
    A = \{\vx\in\bbR^d | x_1\leq -1\},\quad B = \{\vx\in\bbR^d | x_1\geq 1\}.
\end{align}
Now \eqref{eq:DW-potential} is a double-well potential along dimension $x_1$, and the two boundaries $\partial A$ and $\partial B$ are located in the potential wells. When the temperature $T = 1/\beta$ is low, the equilibrium density $p \propto e^{-\beta V}$ is concentrated within the two wells. Meanwhile, in this case $q$ is mostly flat with a sharp transition from $0$ to $1$ at $x_1 = 0$.

For this example, we can compute a ground truth solution. By symmetry, we can obtain the committor function by solving the backward Kolmogorov equation in the first dimension, i.e., setting $q_{\text{true}}(\vx) = f(x_1)$, where
\begin{align}
    \frac{d^2 f(x_1)}{dx_1^2} -4x_1(x_1^2 - 1)\frac{df(x_1)}{dx_1}=0,\quad f(-1)=0,\quad f(1) = 1.
\end{align}
We can solve this ODE numerically using a finite difference method on a very fine grid to produce $q_{\text{true}}$. The performance of our proposed method is evaluated by the relative error metric 
\begin{align}
    E = \frac{\|q-q_{\text{true}}\|_{L^2(\OAB,p)}}{\|q_{\text{true}}\|_{L^2(\OAB,p)}},
    \label{eq:relative-error}
\end{align}
where $\|\cdot\|_{L^2(\OAB,p)}$ denotes the $L^2$-norm with respect to the equilibrium density $p$ defined in \eqref{eq:equi-measure} over the domain $\OAB$. 

We enforce the boundary conditions by constructing soft boundary measures $p_A$ and $p_B$ in MPS/TT format following \eqref{eq:gaussian-for-hyperplane}. Meanwhile, we can exactly treat the equilibrium density $p$ in MPS/TT format since it factorizes as a pure tensor product $$p(\bm{x}) = \prod_{k=1}^d p_k (x_k)$$ of univariate functions, given the choice of potential \eqref{eq:DW-potential}.

It remains to fix a univariate basis for each dimension of the committor function $q$. One could of course choose a generic basis such as Chebyshev polynomials, Legendre polynomials, or Fourier series. For this example, however, a better choice is to construct an appropriate truncated orthogonal polynomial basis for each dimension $k$ according to the univariate density $p_k$.

\begin{figure}[!htb]
    \centering
    \subfloat[]{{\includegraphics[width=0.43\textwidth]{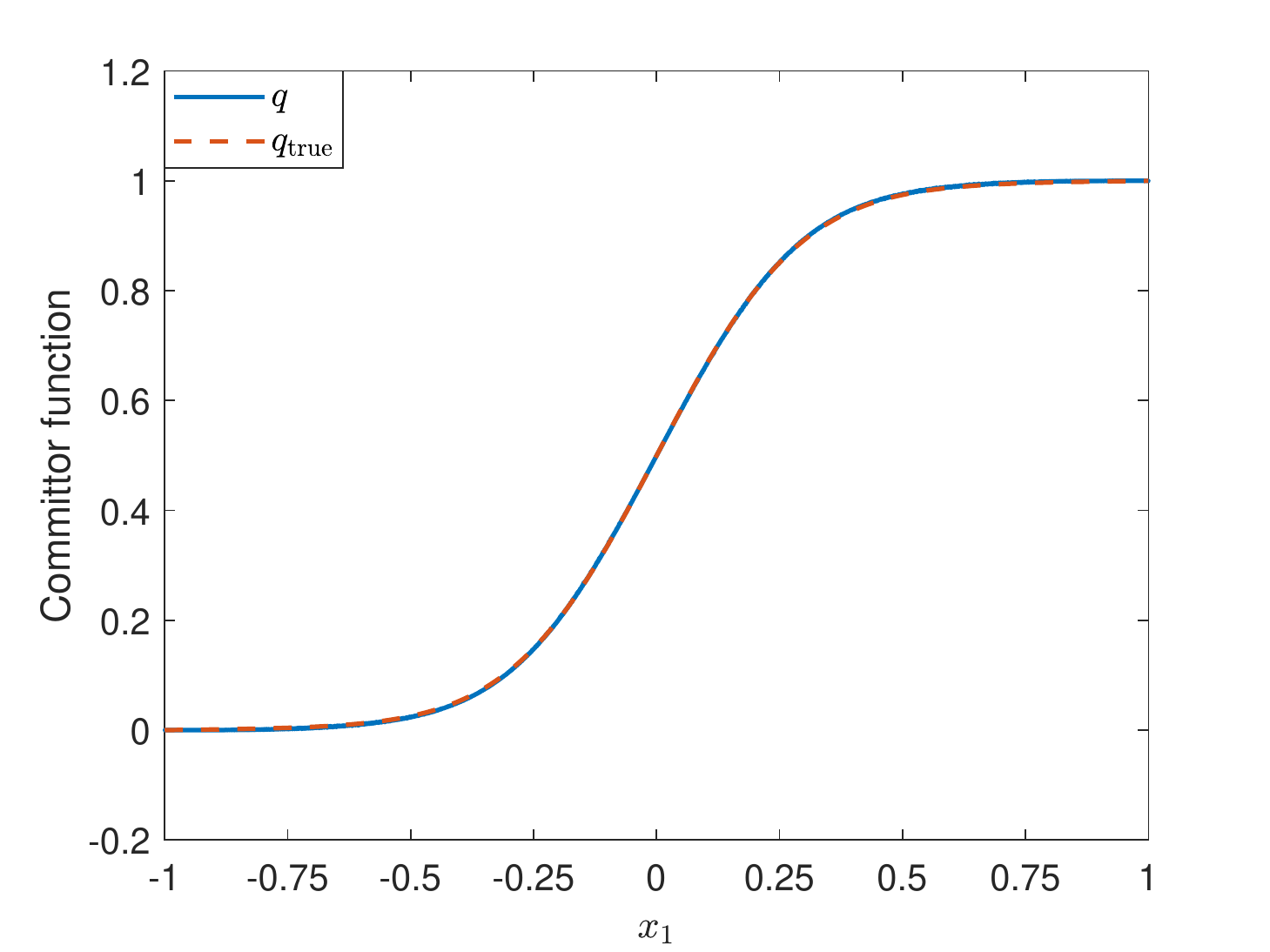}}}\quad
    \subfloat[]{{\includegraphics[width=0.43\textwidth]{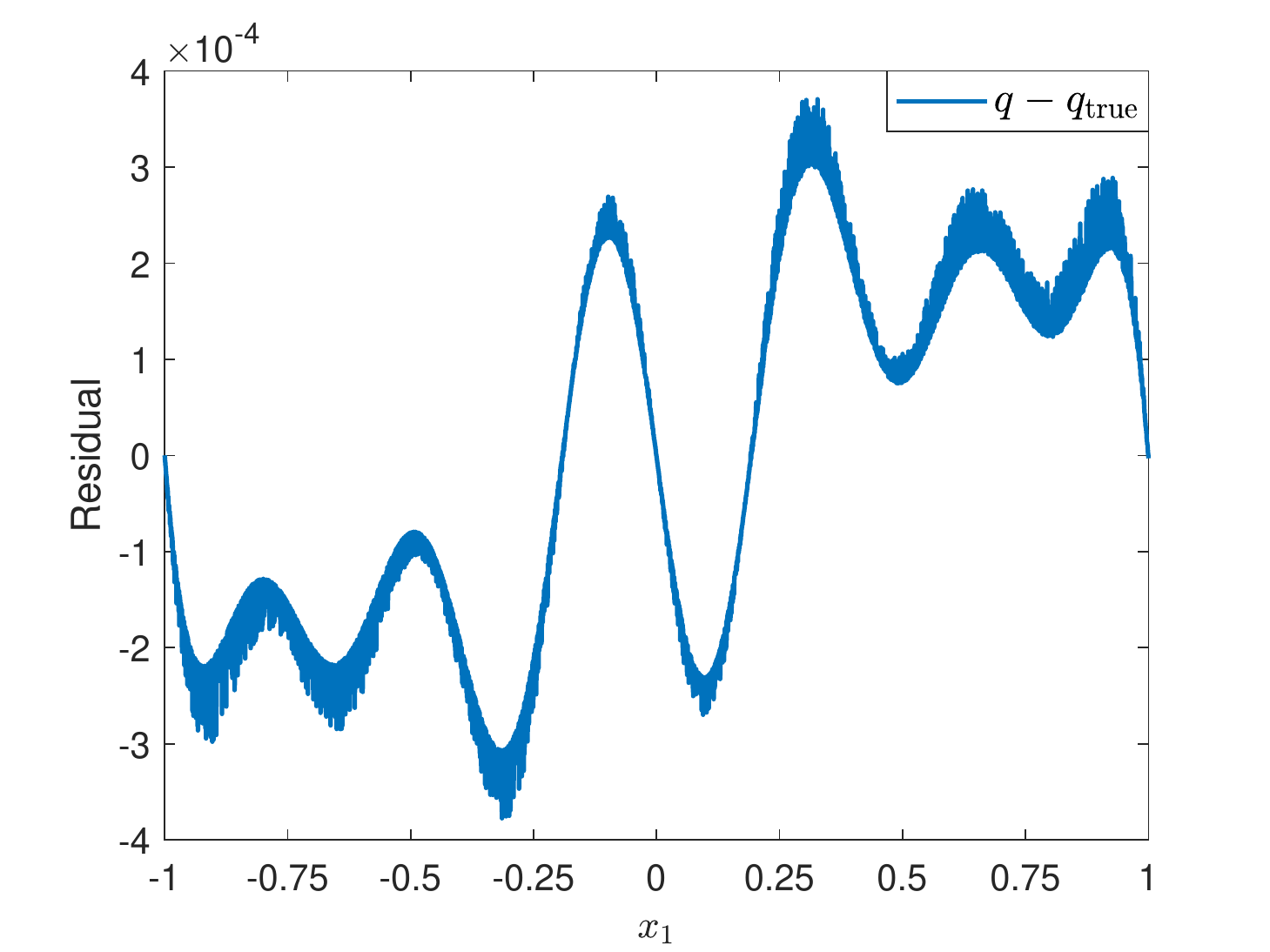}}}
    \caption{Numerical results for the double-well potential, $T=0.2$. (a) The numerical solution of the committor function $q$, compared with the ground truth $q_{\text{true}}$, plotted along the $x_1$ dimension. (b) The residual plot $q-q_{\text{true}}$ along the $x_1$ dimension.}
    \label{fig:DW-solutions}
\end{figure}

\begin{figure}[!htb]
    \centering
    \subfloat[]{{\includegraphics[width=0.43\textwidth]{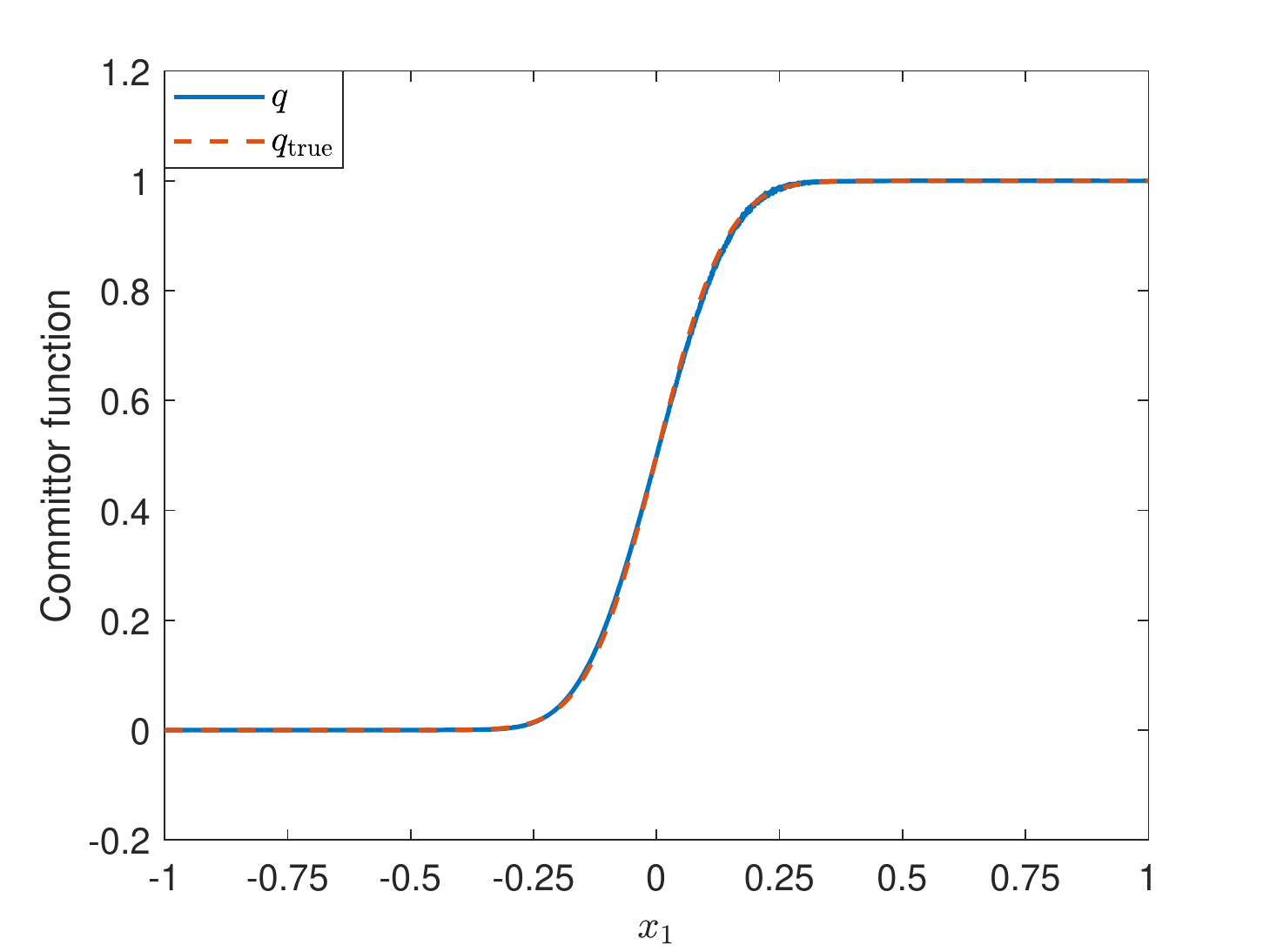}}}\quad
    \subfloat[]{{\includegraphics[width=0.43\textwidth]{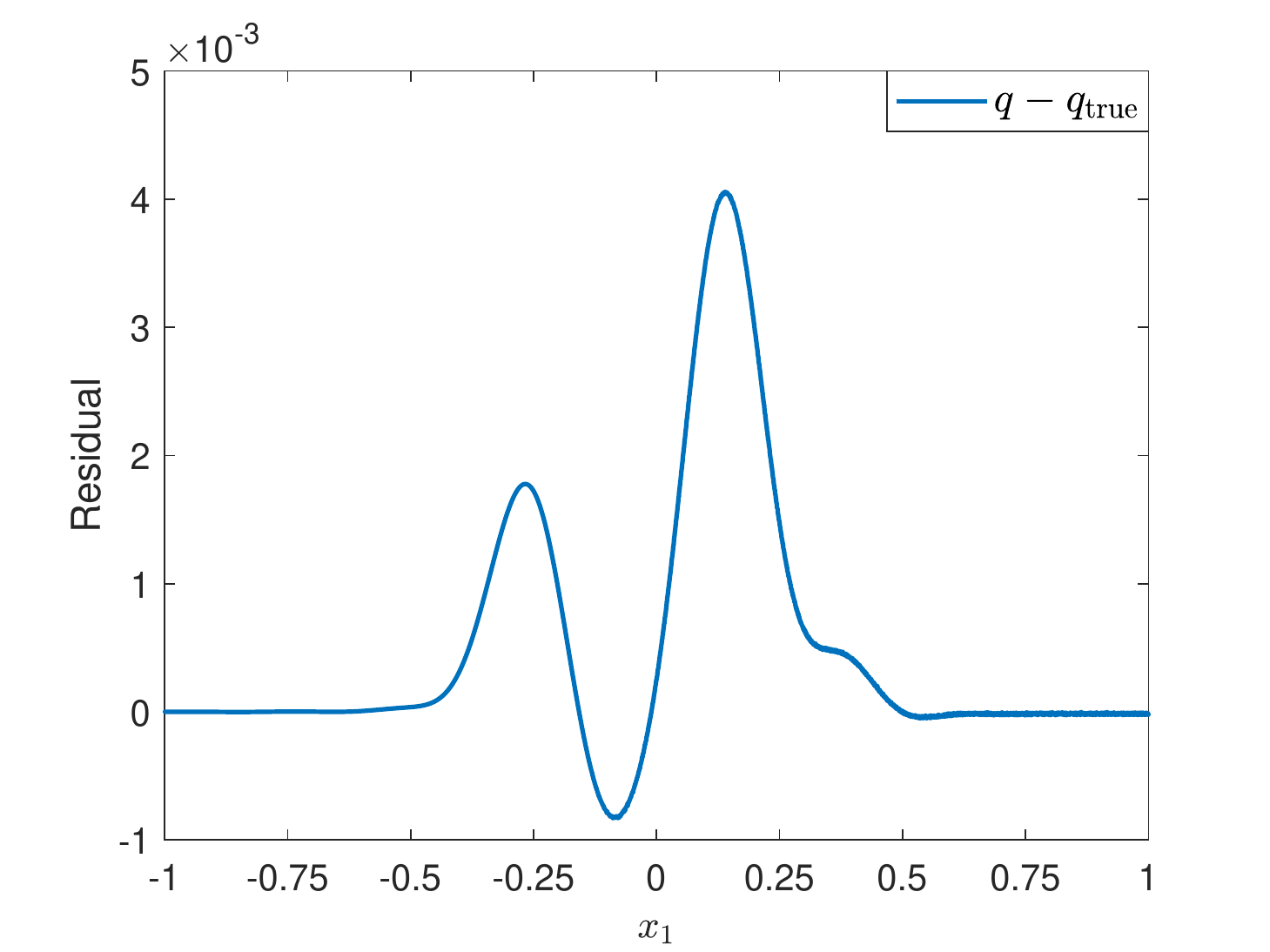}}}
    \caption{Numerical results for the double-well potential, $T=0.05$. (a) The numerical solution of the committor function $q$, compared with the ground truth $q_{\text{true}}$, plotted along the $x_1$ dimension. (b) The residual plot $q-q_{\text{true}}$ along the $x_1$ dimension.}
    \label{fig:DW-solutions2}
\end{figure}

We solve for the committor function at two representative temperatures $T=0.2$ and $T=0.05$ in $d=20$ dimensions. For $T=0.2$, we use the first $30$ orthongonal polynomials for all dimensions. For the lower temperature $T=0.05$, we use $60$ orthogonal polynomial basis functions since the true committor function changes more sharply near $x_1=0$. We show $q$ and $q_{\text{true}}$ for $T=0.2$ in \figref{fig:DW-solutions} (a) and the corresponding residual $q-q_{\text{true}}$ in \figref{fig:DW-solutions} (b). Numerical results for $T=0.05$ are illustrated similarly in \figref{fig:DW-solutions2}.

We compute the relative error \eqref{eq:relative-error} by Monte Carlo integration. In practice, to obtain Monte Carlo error under $10^{-6}$ it suffices to generate $10^7$ samples from the equilibrium density $p$. The relative error is $E=1.60\times 10^{-4}$ for $T=0.2$ and $E=6.77\times 10^{-4}$ for $T=0.05$. Additional tests were performed with other bases such as Chebyshev polynomials and Fourier series. Qualitatively, the behavior was similar, though the performance, as measured via the relative error \eqref{eq:relative-error}, was slightly worse.

\subsection{Ginzburg-Landau potential}

\begin{figure}[!htb]
    \centering
    \subfloat[Local minimizer $U_{-}$]{{\includegraphics[width=0.43\textwidth]{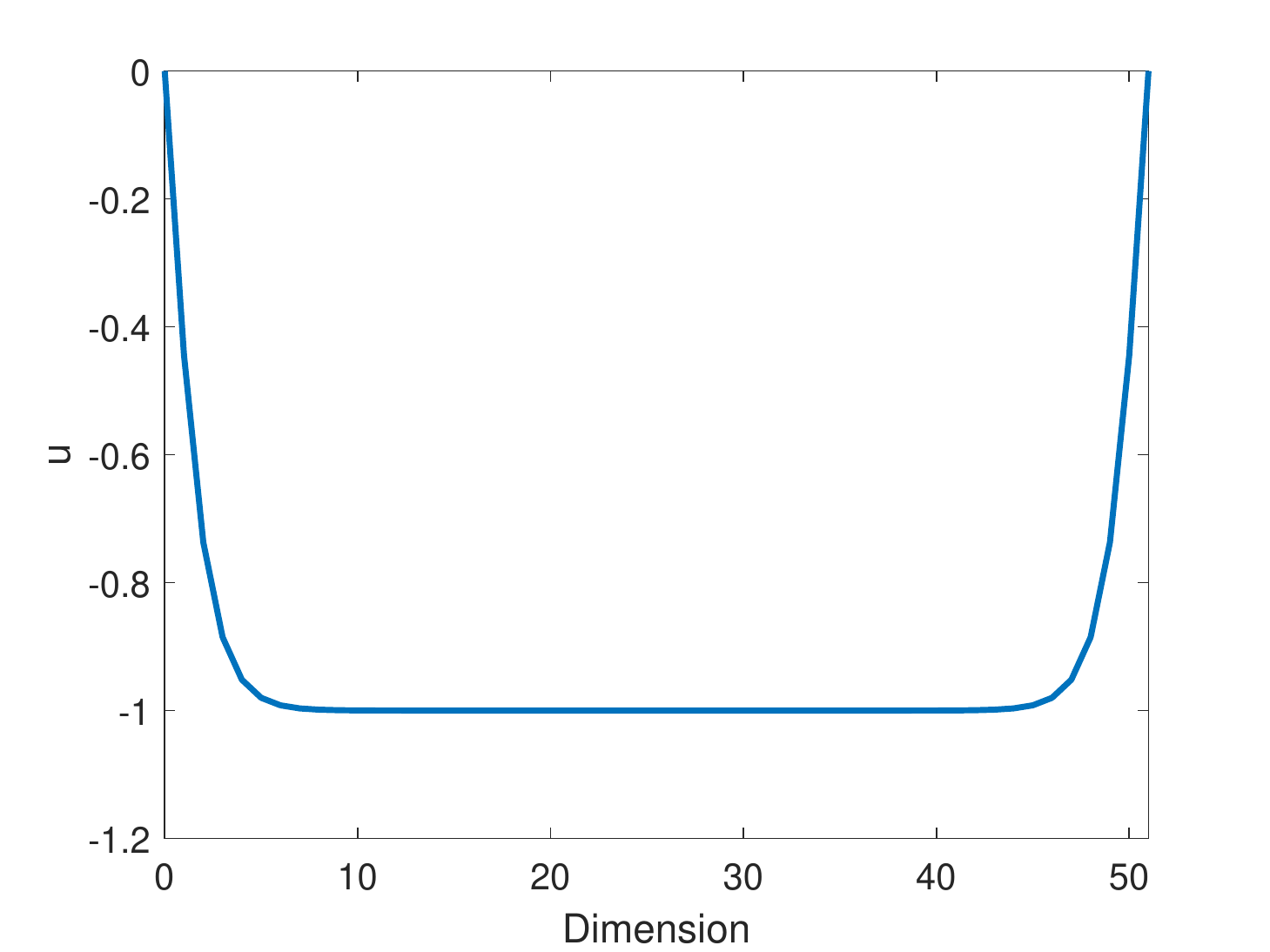}}}\quad
    \subfloat[Local minimizer $U_{+}$]{{\includegraphics[width=0.43\textwidth]{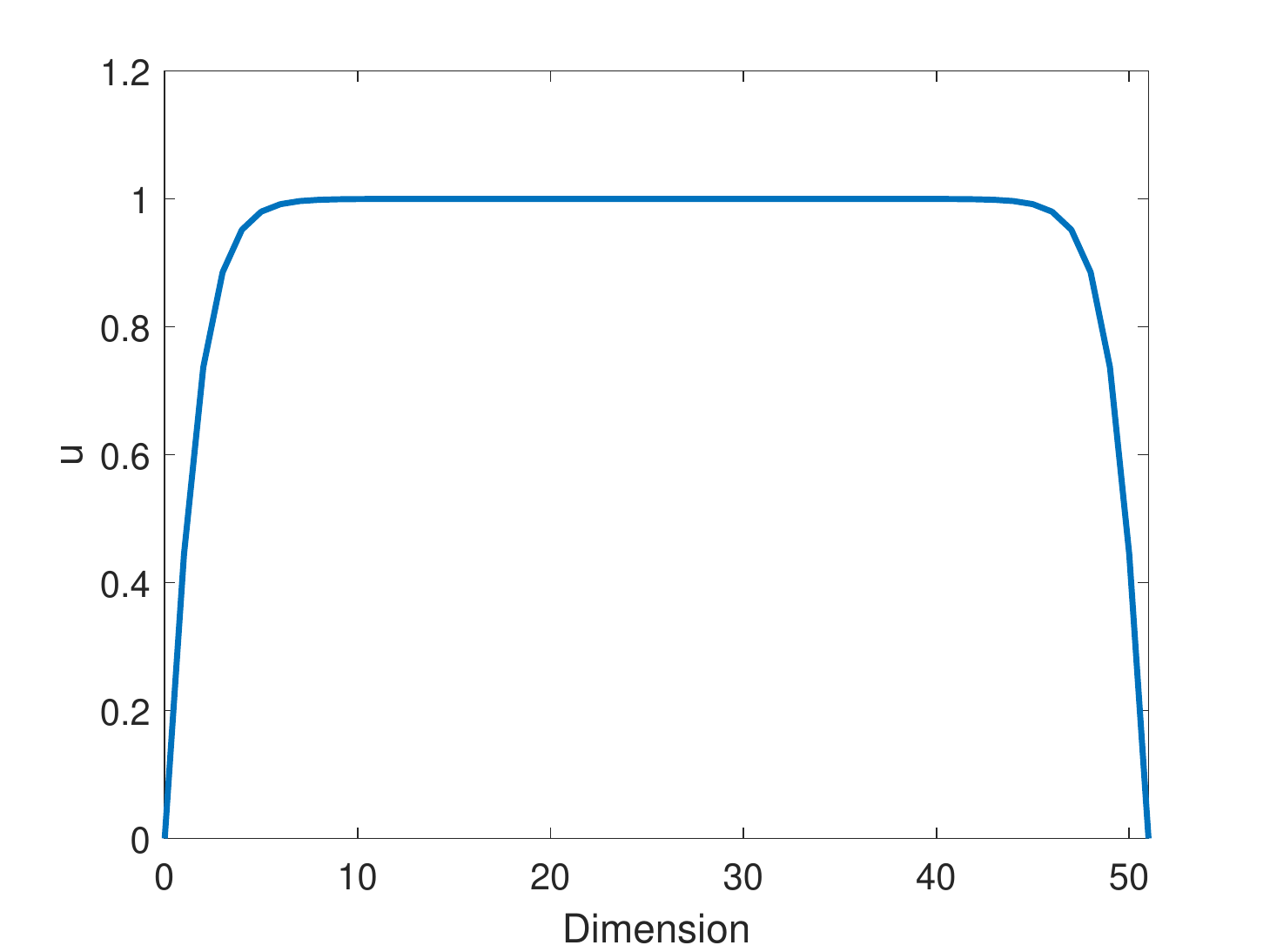}}}
    \caption{The two global minima of the Ginzburg-Landau energy \eqref{eq:discrete-GL} with $d=50$ and $\lambda=0.03$.}
    \label{fig:GL-minima}
\end{figure}

The Ginzburg-Landau theory was developed to provide a mathematical description of superconductivity \cite{GL-model}. In this numerical example, we consider a simplified Ginzburg-Landau model, in which the Ginzburg-Landau energy is defined for a one-dimensional scalar field $u:[0,1]\ra \R$ as follows:
\begin{align}
    \tilde{V} [u] = \int_{0}^1 \left[ \frac{\lambda}{2} (u')^2 + \frac{1}{4\lambda} (1-u^2)^2 \right] dx,
    \label{eq:continuous-GL}
\end{align}
where $\lambda$ is a small positive parameter and $u$ satisfies the boundary conditions $u(0)=u(1)=0$. We discretize $u$ uniformly on $[0,1]$ as $U=(U_1,U_2,\dots,U_d)$ with boundary conditions $U_0=U_{d+1}=0$. Then we approximate the continuous Ginzburg-Landau energy \eqref{eq:continuous-GL} with the discretization 
\begin{align}
    V(U) := \sum_{i=1}^{d+1} \frac{\lambda}{2} \left(\frac{U_i - U_{i-1}}{h}\right)^2 + \frac{1}{4\lambda} (1 - U_i^2)^2,
    \label{eq:discrete-GL}
\end{align}
where the grid spacing $h = 1/(d+1)$. We fix $d=50$ and $\lambda=0.03$. Note that $V(U)$ has two global minima $U_{\pm}$ illustrated in \figref{fig:GL-minima}. We let $A$ and $B$ be the balls
$\{U:\|U- U_{\pm}\| \leq R \}$
centered at the global minima. The radius $R$ is set to be $2.5$, chosen such that the balls $A$ and $B$ roughly contain the regions of high equilibrium probability density around the two centers.

We present numerical results for two representative temperatures $T=8$ and $T=16$. We enforce the boundary conditions by constructing soft boundary measures $p_A$ and $p_B$ in MPS/TT format following \eqref{eq:high-d-gaussian-for-sphere}. We detail the approximation of the equilibrium probability density $p$ in MPS/TT format in \appref{app:GL-density}. 

We define the domain to be the hypercube $\Omega=[-\gamma,\gamma]^{50}$. Based on our choices of $\lambda$ and $\beta$, we take $\gamma=2.6$ since this choice guarantees that the equilibrium density has negligible mass outside of $\Omega$. To represent the committor function $q$, we use the first $5$ Fourier basis functions $\{1,\cos{(\pi x/\gamma)},\sin{(\pi x/\gamma)},\cos{(2\pi x/\gamma)},\sin{(2\pi x/\gamma)}\}$ for each dimension. The ranks of the coefficient MPS/TT $\ccQ$ are all taken to be $6$. We initialize all the entries of the unknown tensor cores of $\ccQ$ with normal $\mathcal{N}(0,1)$ random numbers and then perform ALS, gradually increasing the penalty parameter $\rho$ to better enforce the boundary conditions. In practice we observe that different initializations have little effect on the output of the algorithm.

For problems of this size, traditional methods are intractable, making it difficult to obtain an exact reference $q_{\text{true}}$ for comparison. Instead, as a proxy we study a `thickened isosurface' around $q = 0.5$, defined as $\Gamma_{\epsilon} = \{U:\|q(U) - 0.5\|\leq \epsilon\}$, where $\epsilon > 0$ is a small threshold parameter. If the solution $q$ is indeed a satisfactory approximation of the true committor function, then for any trajectories given by \eqref{eq:langevin} starting from points in $\Gamma_{\epsilon}$, the probability of entering region $B$ before $A$ should be close to $0.5$.

\begin{figure}[!htb]
    \centering
    \subfloat[Histogram]{{\includegraphics[width=0.4\textwidth]{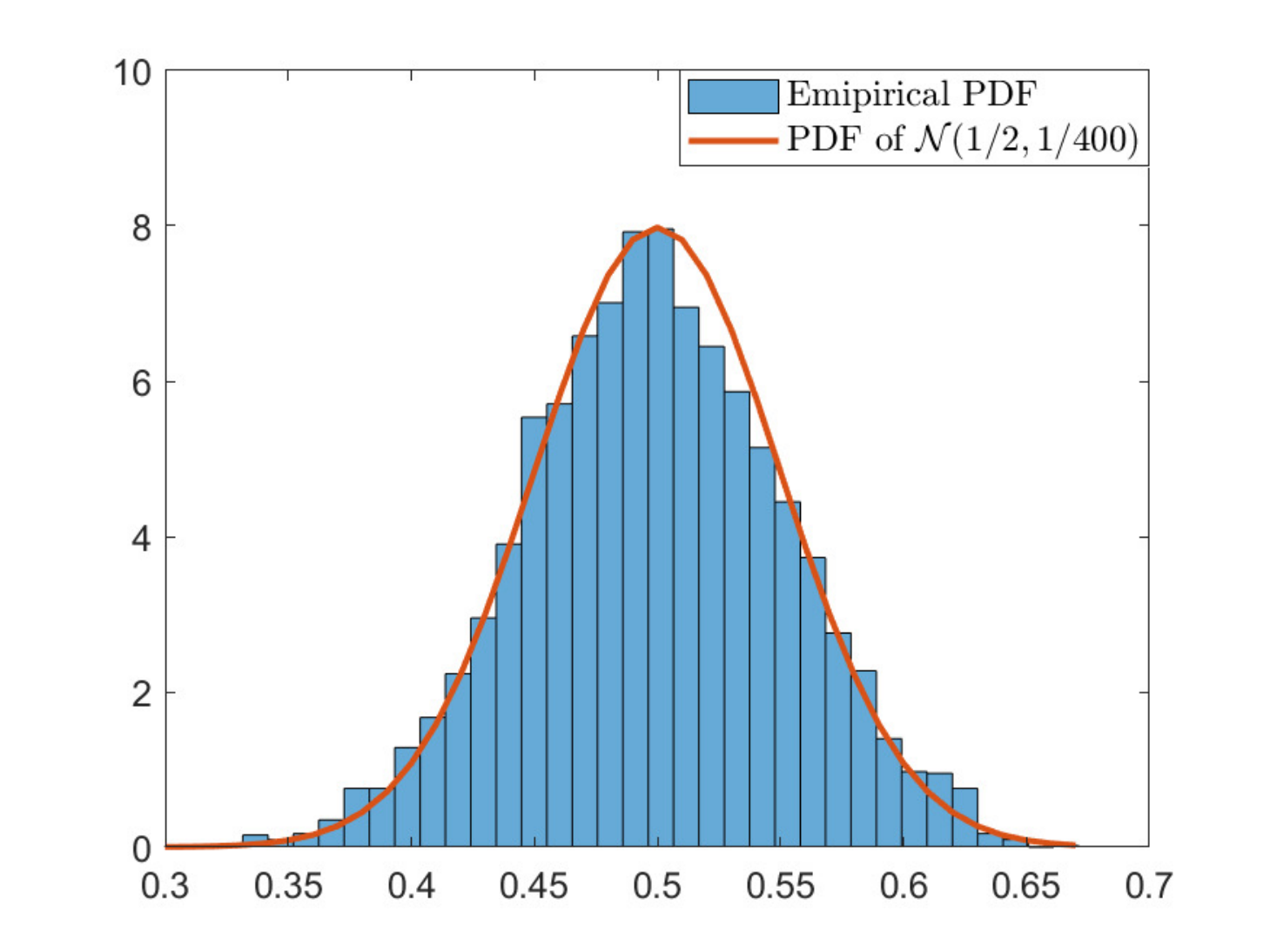}}}
    \subfloat[Q-Q plot]{{\includegraphics[width=0.4\textwidth]{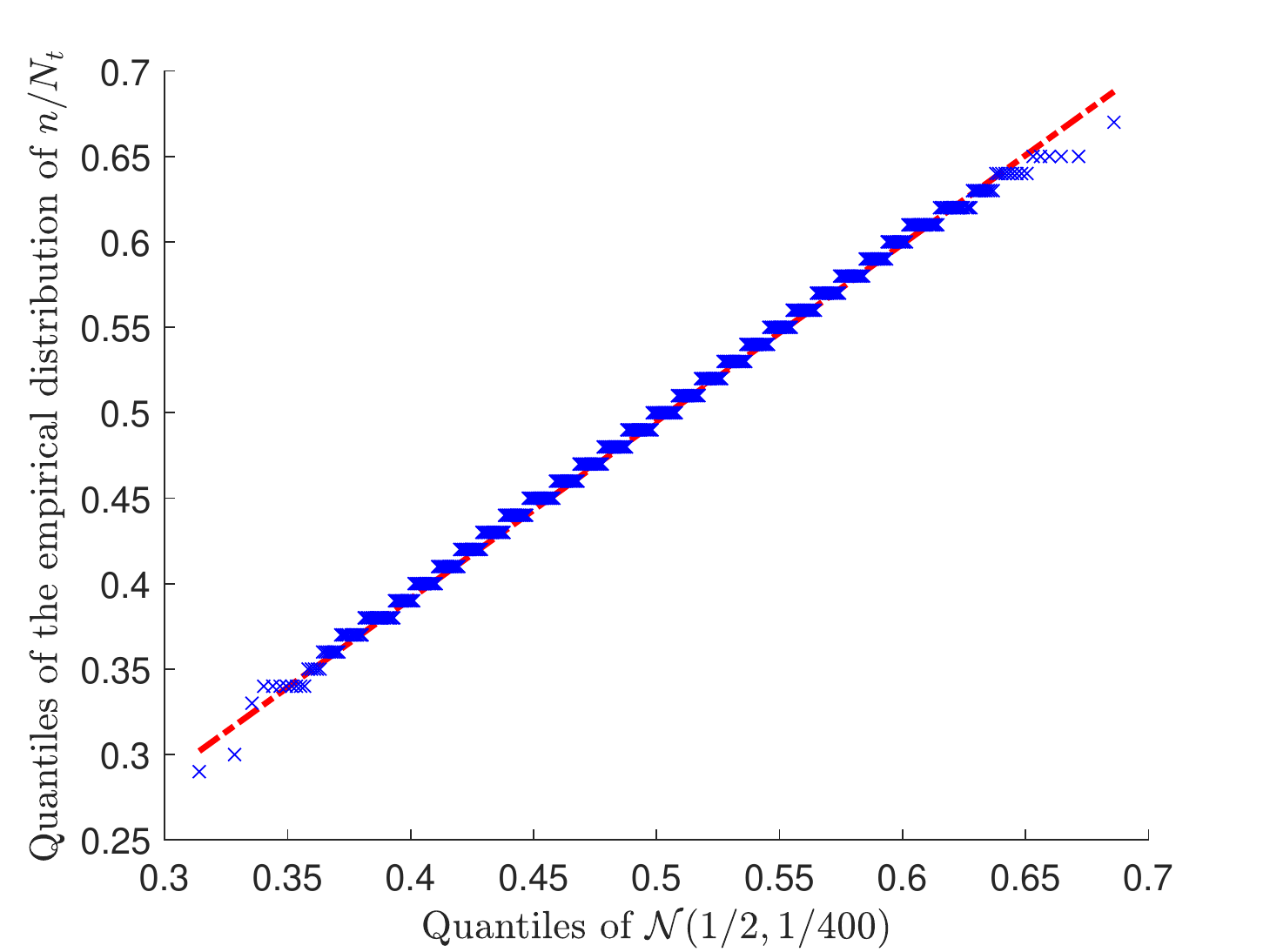}}}
    \caption{Comparison of distributions for $T=8$. (a) Empirical histogram for $\{n_j/N_t\}_{j=1}^{5000}$ compared with the density of $\mathcal{N}(\frac{1}{2},1/400)$. (b) Q-Q (quantile-quantile) plot of $\{n_j/N_t\}_{j=1}^{5000}$ compared with $\mathcal{N}(\frac{1}{2},1/400)$.}
    \label{fig:traj-T8}
\end{figure}

\begin{figure}[!htb]
    \centering
    \subfloat[Probability density]{{\includegraphics[width=0.4\textwidth]{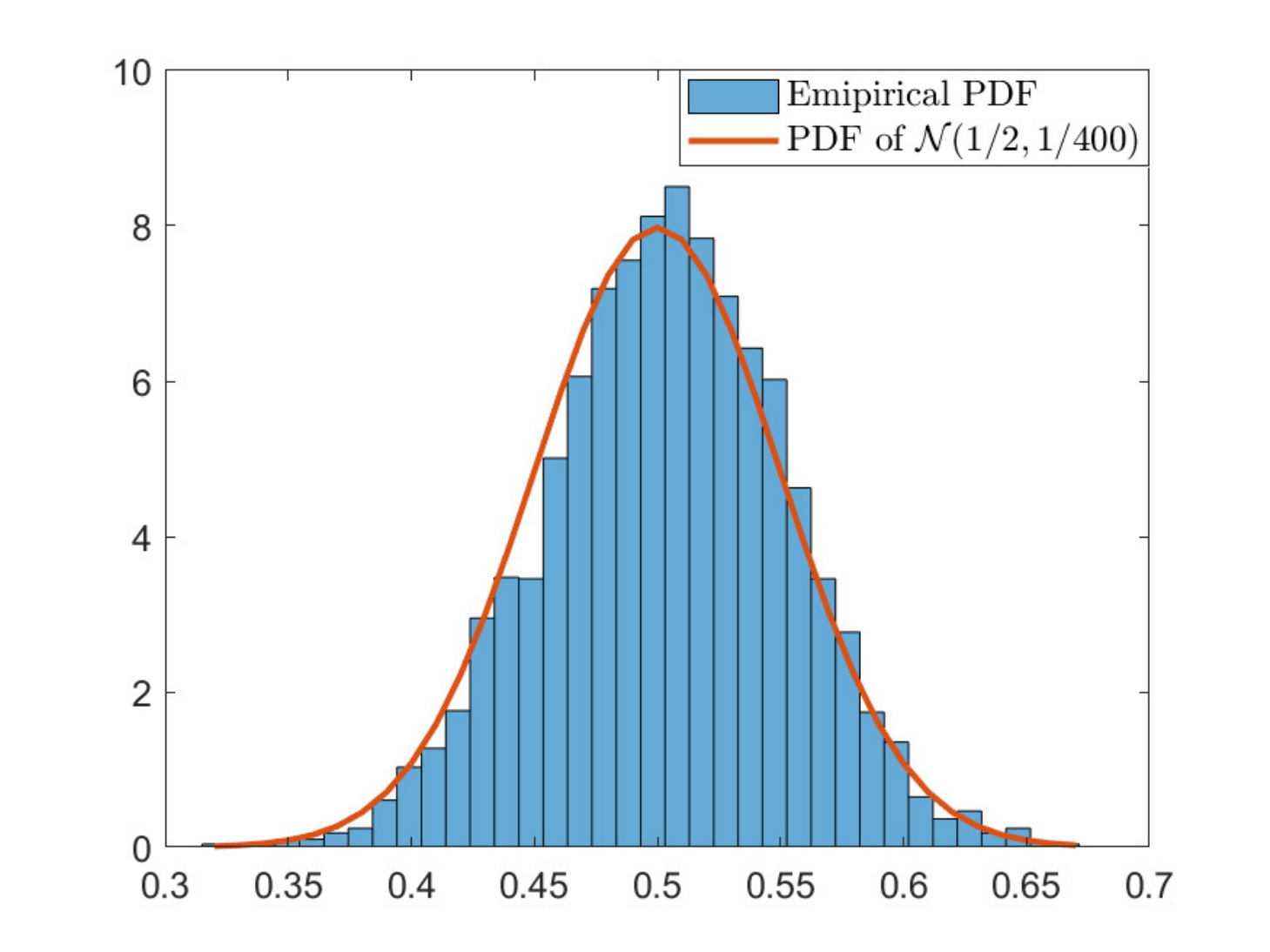}}}
    \subfloat[Q-Q plot]{{\includegraphics[width=0.4\textwidth]{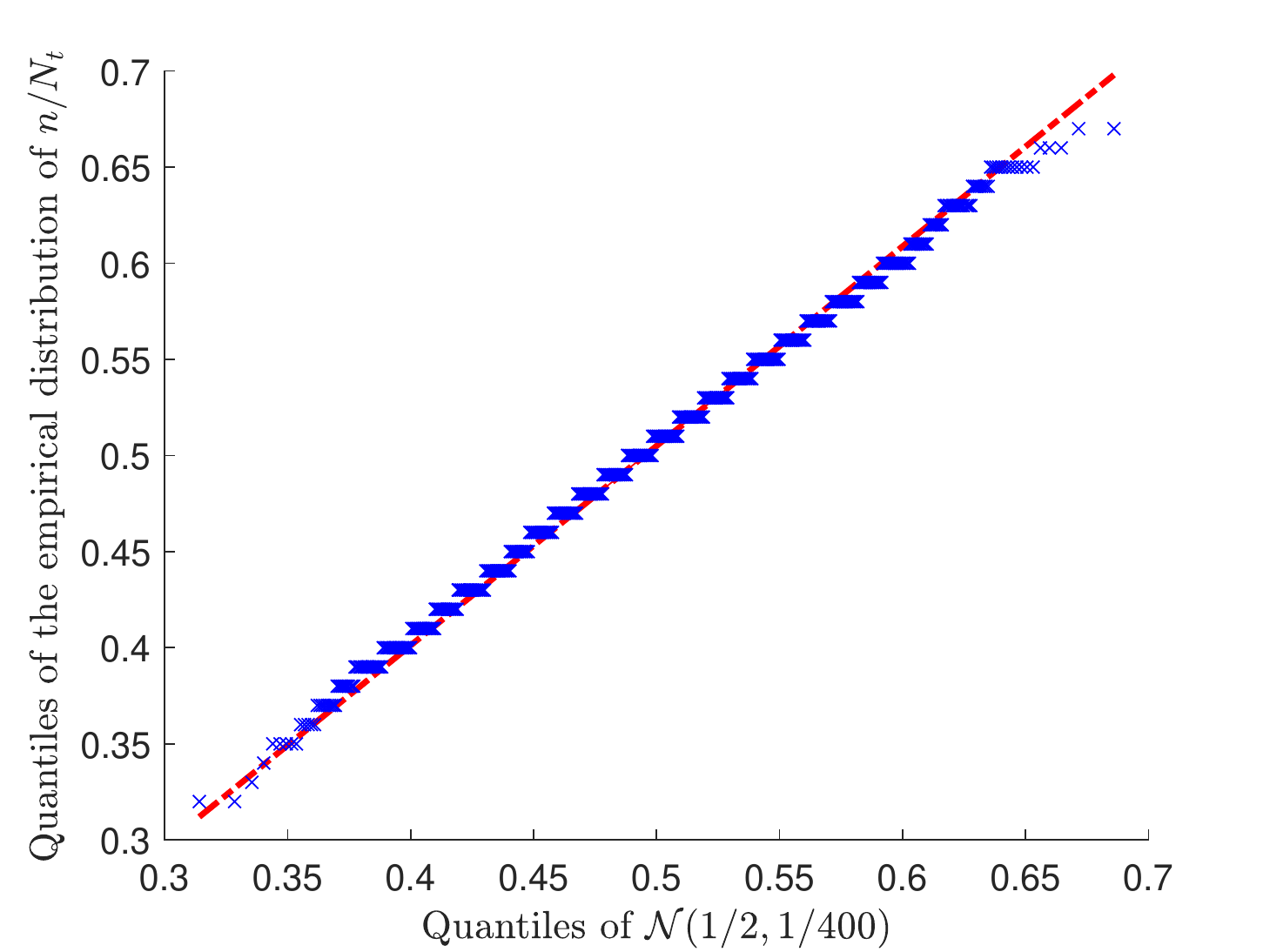}}}
    \caption{Comparison of distributions for $T=16$. (a) Empirical histogram for $\{n_j/N_t\}_{j=1}^{5000}$ compared with the density of $\mathcal{N}(\frac{1}{2},1/400)$. (b) Q-Q (quantile-quantile) plot of $\{n_j/N_t\}_{j=1}^{5000}$ compared with $\mathcal{N}(\frac{1}{2},1/400)$.}
    \label{fig:traj-T16}
\end{figure}

To verify this, we generate samples from the equilibrium distribution by simulating the process \eqref{eq:langevin}. Then we filter to keep samples on the isosurface $\Gamma_\epsilon$ using the computed committor function $q$ and the threshold $\epsilon$ of our choice. Let us pick $N_s$ points in $\Gamma_{\epsilon}$, denoted $\{\tilde{U}_j\}_{j=1}^{N_s}$. For each point $\tilde{U}_j$, we generate $N_t$ trajectories by simulating the Langevin process \eqref{eq:langevin} and use $n_j$ to denote the number of trajectories ending up in region $B$ before $A$. By the central limit theorem, when $N_t$ is large, the distribution of $n_j/N_t$ should be well-approximated by the normal distribution $\mathcal{N}(\frac{1}{2},(4N_t)^{-1})$. In our numerical tests, we set $\epsilon=5\times 10^{-3}$, $N_s=5000$, and $N_t=100$. The results for $T=8$ and $T=16$ are illustrated in \figref{fig:traj-T8} and \figref{fig:traj-T16}, respectively. We compare the histogram of $\{n_j/N_t\}_{j=1}^{5000}$ with the normal distribution $\mathcal{N}(\frac{1}{2},1/400)$ on the left and show the Q–Q (quantile-quantile) plot of the distribution of $\{n_j/N_t\}_{j=1}^{5000}$ versus $\mathcal{N}(\frac{1}{2},1/400)$ on the right. These figures demonstrate that the distribution of $\{n_j/N_t\}_{j=1}^{5000}$ is indeed in good agreement with the normal distribution $\mathcal{N}(\frac{1}{2},1/400)$, which indicates that our solution $q$ provides a good approximation of the true isosurface.

\begin{figure}[!htb]
    \centering
    \subfloat[Two centroids]{{\includegraphics[width=0.4\textwidth]{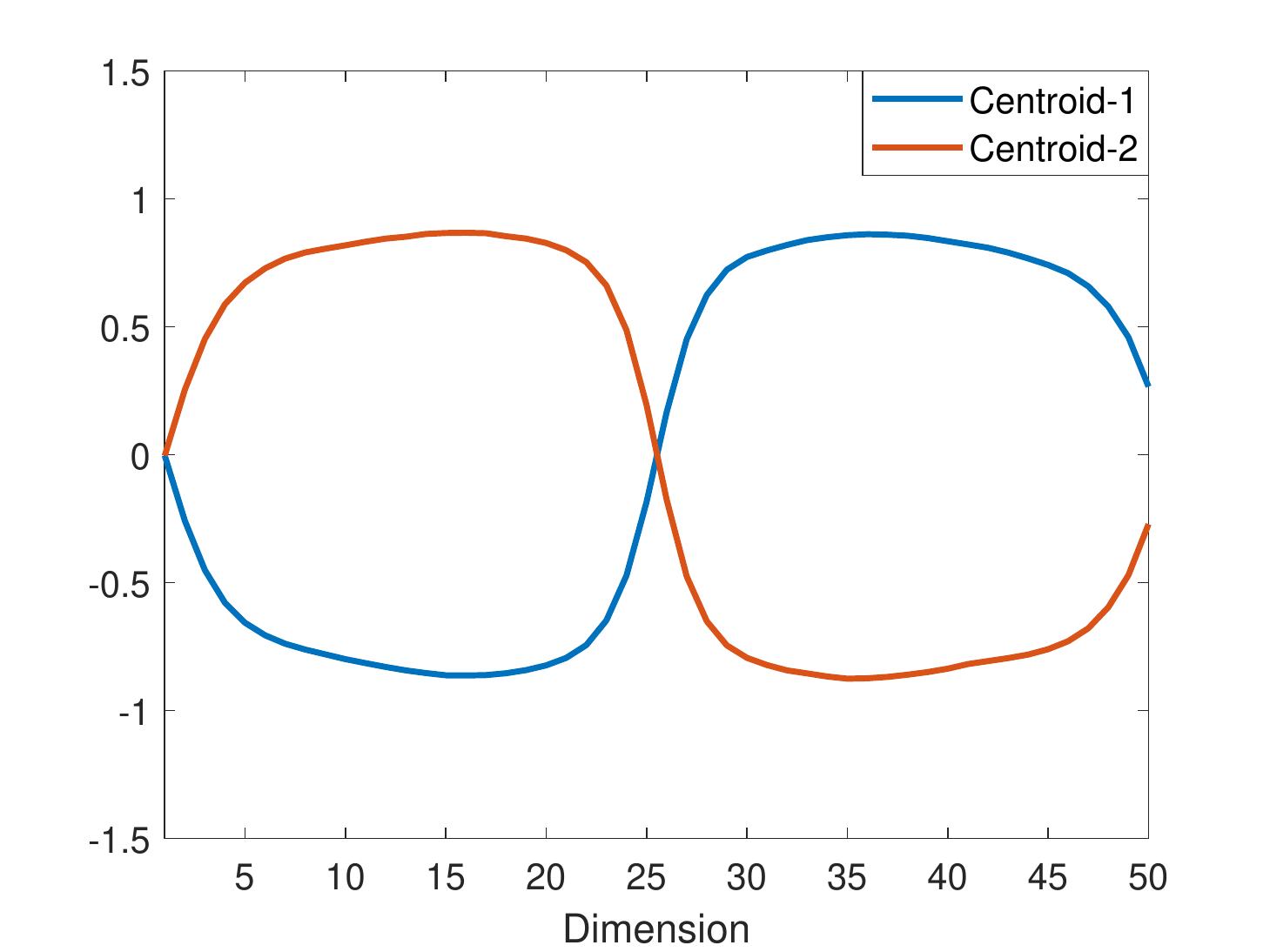}}}
    \subfloat[Histogram of $\theta$]{{\includegraphics[width=0.4\textwidth]{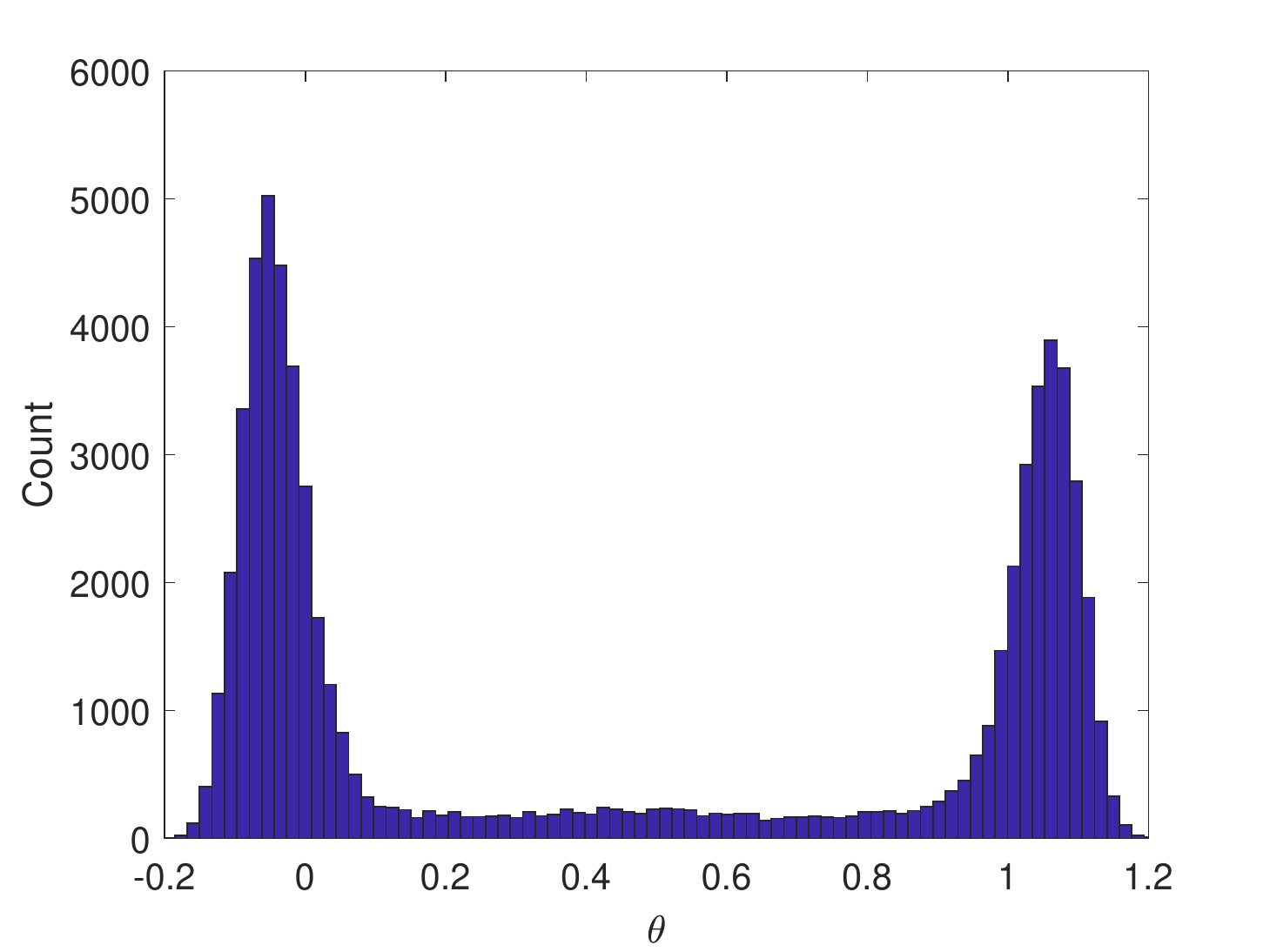}}}
    \caption{Analysis of the $q=0.5$ isosurface for $T=8$. (a) Centroids of the two clusters in the $q=0.5$ isosurface. (b) Histogram of the $1$-dimensional coordinate $\theta$ of the isosurface samples along the line between the two clusters.}
    \label{fig:T8}
\end{figure}

\begin{figure}[!htb]
    \centering
    \subfloat[Two centroids]{{\includegraphics[width=0.4\textwidth]{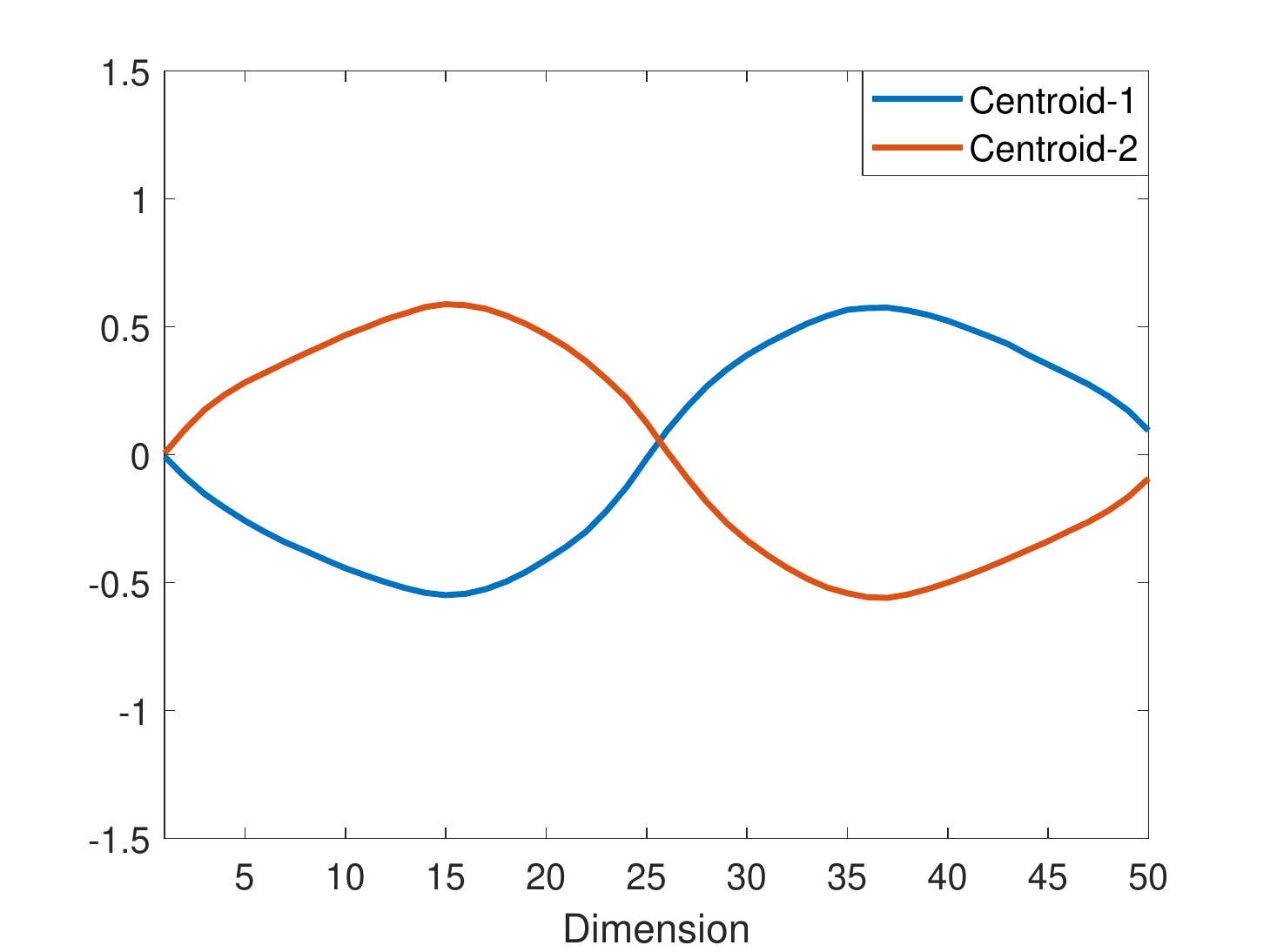}}}
    \subfloat[Histogram of $\theta$]{{\includegraphics[width=0.4\textwidth]{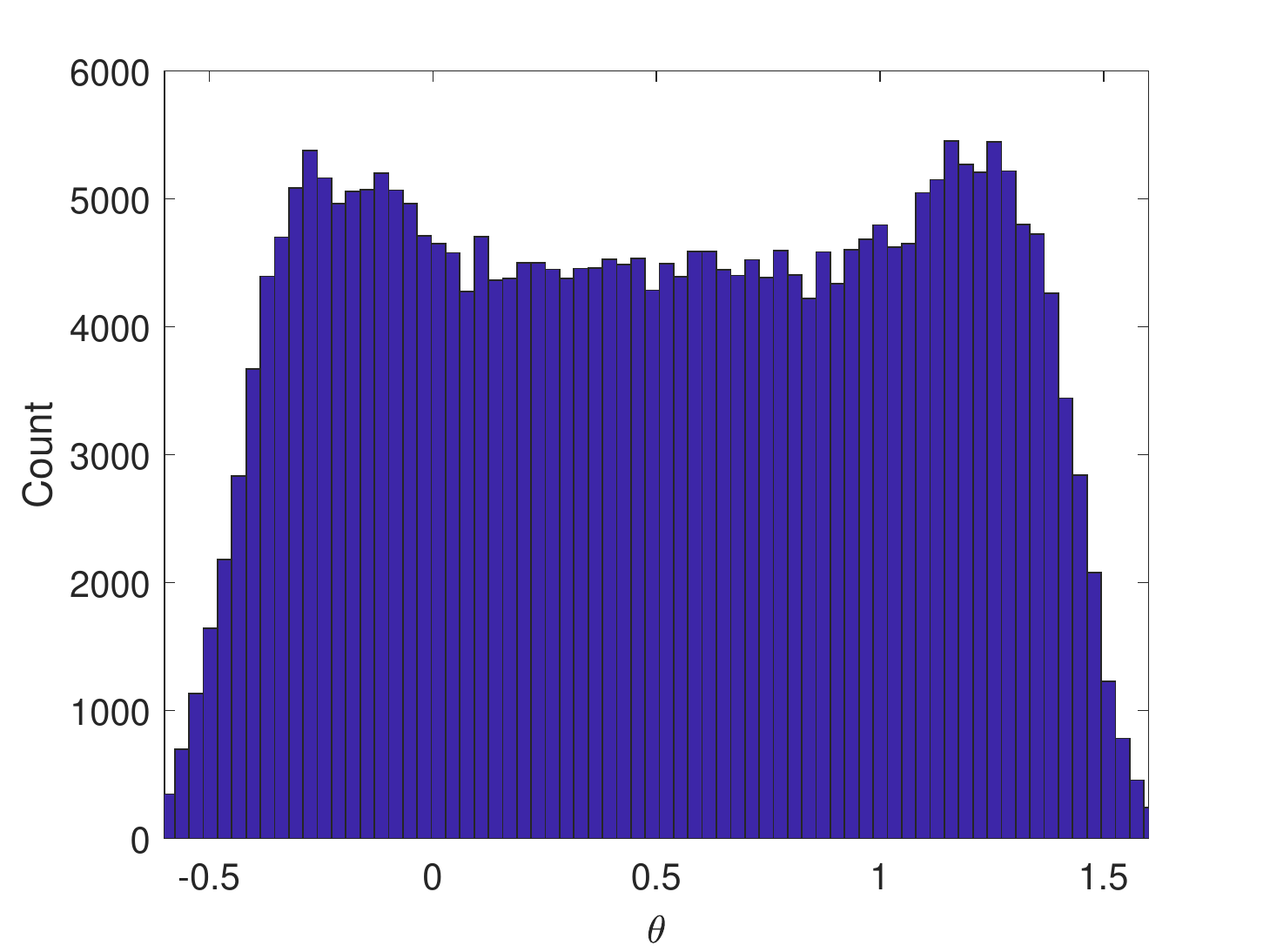}}}
    \caption{Analysis of the $q=0.5$ isosurface for $T=16$. (a) Centroids of the two clusters in the $q=0.5$ isosurface. (b) Histogram of the $1$-dimensional coordinate $\theta$ of the isosurface samples along the line between the two clusters.}
    \label{fig:T16}
\end{figure}

Next, consider the restriction of the equilibrium density to the $q=0.5$ isosurface. Intuitively, the first term in the Ginzburg-Landau potential \eqref{eq:discrete-GL} encourages the configuration $U$ to be as flat as possible. Therefore, to transition between the two boundary states $U_{-}$ and $U_{+}$, it is favorable in terms of energy to have only a single sign change in the discretized function $U$. We compute $10^7$ samples from the equilibrium density by running the overdamped Langevin process \eqref{eq:langevin}, initialized at random states in $\Omega$. We retain only the samples that fall in the thickened isosurface $\Gamma_{\epsilon}$,  $\epsilon=5\times 10^{-3}$. Then we perform $2$-means clustering on these samples. In \figref{fig:T8} (a) and \figref{fig:T16} (a), we show the centroids $U^{(1)}$ and $U^{(2)}$ of the two clusters for $T=8$ and $T=16$, respectively. These configurations are symmetric with a single sign change.

Next we project all samples in the $q=0.5$ isosurface to the line containing the two centroids, i.e., to points of the form $\theta U^{(1)}+(1-\theta)U^{(2)}$. In \figref{fig:T8} (b) and \figref{fig:T16} (b) we plot the histograms of $\theta$ for all samples to demonstrate that these distributions are indeed bimodal. Observe that at higher temperature, the bimodality is less pronounced.

\begin{figure}[!htb]
    \centering
    \subfloat{{\includegraphics[width=0.5\textwidth]{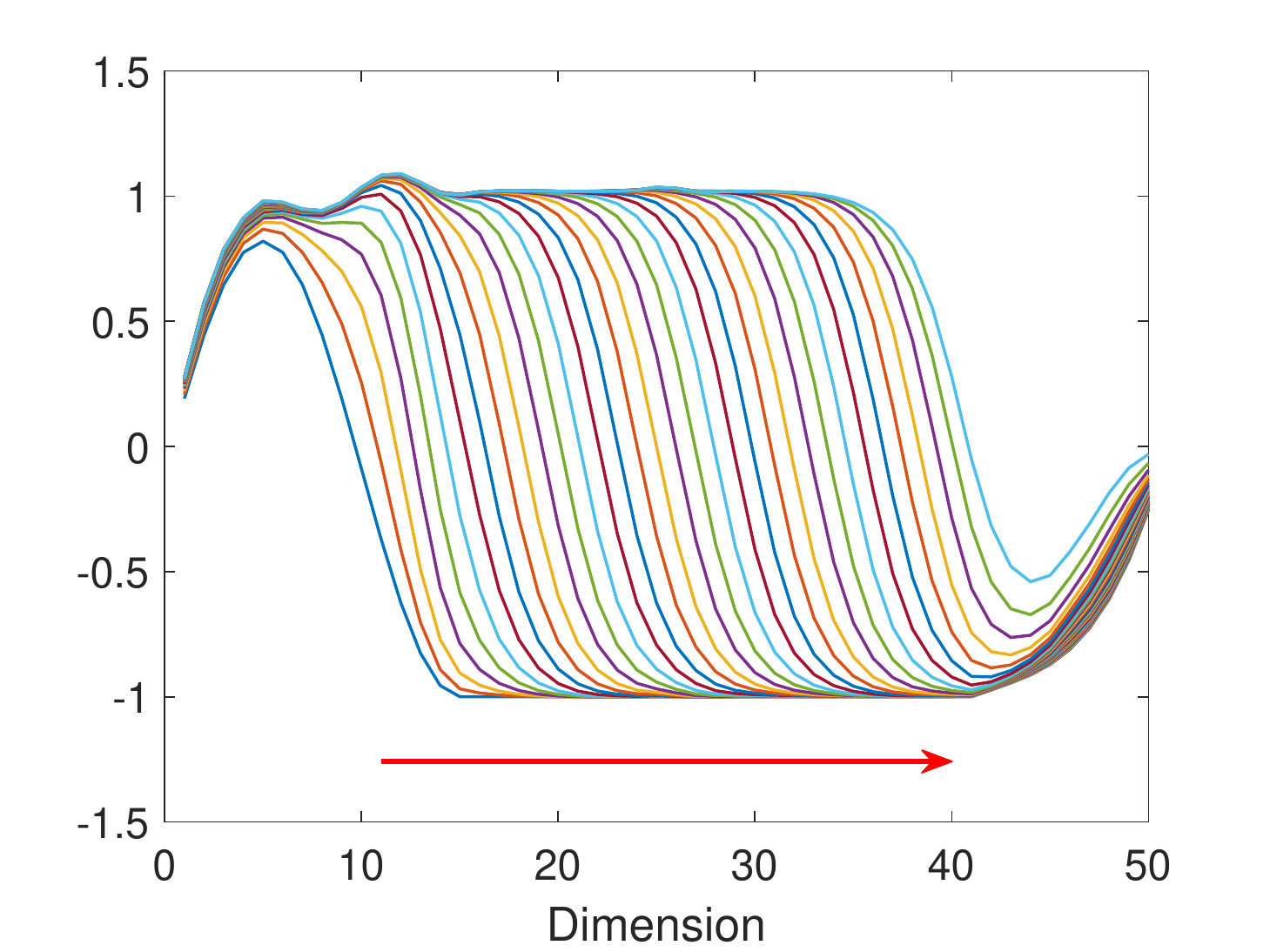}}}
    \caption{Visualization of one transition path for $T=8$.}
    \label{fig:transition-T8}
\end{figure}

Finally, we study transition paths via the deterministic reactive flow \cite{point-committor}:
\begin{align}
    \frac{dU(t)}{dt} = \frac{1}{\beta} p(U(t)) \nabla q(U(t)).
    \label{eq:transition-ODE}
\end{align}
Based on \figref{fig:T8}, we expect that at low temperatures the transition paths between $A$ and $B$ are localized within in two reaction tubes. We visualize one of the transition paths at temperature $T=8$ in \figref{fig:transition-T8}. The leftmost curve corresponds to the initial state of  \eqref{eq:transition-ODE}, for which  $q=0.1$. Meanwhile $q=0.9$ for the rightmost curve. The red arrow indicates the direction of time evolution.

\section{Conclusion}
\label{sec:conclusion}

In this paper, we propose a novel approach for computing high-dimensional committor functions using MPS/TT. In particular, we start from the variational formulation \eqref{eq:soft-varational} for the soft committor function, which can be viewed as an 
approximation of the committor function but which also enjoys a probabilistic interpretation in its own right. To compute high-dimensional integrals, we approximate the equilibrium density and soft boundary measures in MPS/TT format. Meanwhile, the unknown committor function $q$ is also parametrized in MPS/TT format. The variational problem can then be reformulated using standard MPO and MPS/TT operations, and the optimization of $q$ can be performed with $\ccO(d)$ complexity. Extensive numerical experiments demonstrate the computational efficiency and accuracy of the proposed method.

\newpage
\appendix

\section{Soft committor functions}
\label{app:soft}
Consider the optimization problem \eqref{eq:soft-varational}.
We offer a probabilistic interpretation of the optimizer $q$, which we call a `soft committor function.' This probabilistic interpretation will justify its use
qualitatively and quantitatively as a proxy for the committor function.

Define 
\[
 f:= \frac{\rho \cdot p_A}{\beta \cdot p}, \quad  g:= \frac{\rho \cdot p_B}{\beta \cdot  p}, 
\]
and note that the Euler-Lagrange equation associated to (\ref{eq:soft-varational}) then reads as 
\begin{equation}
-\beta^{-1}\Delta q(\bm{x})+\nabla V(\bm{x})\cdot\nabla q(\bm{x})+\left(f(\bm{x})+g(\bm{x})\right)q(\bm{x})=g(\bm{x}).\label{eq:softEL}
\end{equation}
 We shall now rederive this PDE (\ref{eq:softEL}) via a probabilistic
construction. In particular this construction will imply that the
solution $q$ satisfies $0\leq q(\bm{x})\leq1$ for all $\bm{x}\in\Omega$.

Consider a stochastic process $\bm{X}_{t}$ that modifies the standard
overdamped Langevin diffusion 
\[
d\bm{X}_{t}=-\nabla V(\bm{X}_{t})\,dt+\sqrt{2\beta^{-1}}\,d\bm{W}_{t},
\]
 where $\bm{W}_{t}$ is a Wiener processs, by adding jumps to one
of \emph{two} possible `cemetery states' $c_{A},c_{B}$ with state-dependent
rates specified by $f,g$, respectively.

In other words, we view $\bm{X}_{t}\in\Omega\cup\{c_{A}\}\cup\{c_{B}\}$
as the continuous-time limit of the discrete-time Markov chain (also
denoted $\bm{X}_{t}$, abusing notation slightly) defined by the update
\[
\bm{X}_{t+\Delta t}=\begin{cases}
c_{A}, & \mbox{if }U_{t}^{A}<f(\bm{X}_{t})\,\Delta t,\\
c_{B}, & \mbox{otherwise if }U_{t}^{B}<g(\bm{X}_{t})\,\Delta t\\
\bm{X}_{t}-\nabla V(\bm{X}_{t})\,\Delta t+\sqrt{2\beta^{-1}\Delta t}\,\bm{Z}_{t} & \mbox{otherwise},
\end{cases},
\]
 whenever $\bm{X}_{t}\in\Omega$. Here the $U_{t}^{A},U_{t}^{B}$
are i.i.d. uniformly distributed random variables on $[0,1]$, and
the $\bm{Z}_{t}$ are i.i.d. standard Gaussian random variables. In
other words, at each time step the stochastic process is sent to the
cemetery state $c_{A}$ with probability $f(\mathbf{X}_{t})\,\Delta t$
and the cemetery state $c_{B}$ with probability $g(\mathbf{X}_{t})\,\Delta t$,
else it is advanced by the usual overdamped Langevin dynamics. Note
that in the limit of $\Delta t$ small, it is unlikely for both $U_{t}^{A}<f(\bm{X}_{t})\,\Delta t$
and $U_{t}^{B}<g(\bm{X}_{t})\,\Delta t$ to hold. Specifically, the
probability of this is only $O(\Delta t^{2})$ and does not affect
the continuous-time limit.

Moreover, if $\bm{X}_{t}=c_{A}$ (resp., $c_{B})$, then we define
$\bm{X}_{t+\Delta t}=c_{A}$ (resp., $c_{B})$ deterministically.
Then define the stopping time $\tau$ as 
\[
\tau=\inf\left\{ t\,:\,\bm{X}_{t}\in\{c_{A},c_{B}\}\right\} ,
\]
 and define $q:\Omega\ra[0,1]$ by 
\[
q(\bm{x})=\P\left(\bm{X}_{\tau}=c_{B}\,\vert\,\bm{X}_{0}=\bm{x}\right).
\]
 We claim that $q$ so defined (in the continuous-time limit $\Delta t\ra0$)
satisfies the PDE (\ref{eq:softEL}). Note that this construction
of $q$ coincides with the usual probabilistic construction for the
committor function, modulo a change in the underlying stochastic process.
We justify the claim only formally. A rigorous argument can be made
by analogy to arguments made for the ordinary committor function~\cite{bhattacharya2009stochastic}.

\subsection{Formal PDE derivation}
The PDE is derived by conditioning on $\bm{X}_{0}=\bm{x}$, and writing
\begin{eqnarray*}
q(\bm{x}) & = & \P(\bm{X}_{\tau}=c_{B})\\
 & = & \P(\bm{X}_{\Delta t}=c_{B})+\P(\bm{X}_{\Delta t}\notin\{c_{A},c_{B}\})\,\E\left[\P\left(\bm{X}_{\tau}'=c_{B}\,\vert\,\bm{X}_{0}'=\bm{X}_{\Delta t}\right)\right]+O(\Delta t^{2}),
\end{eqnarray*}
 where $\bm{X}_{t}'$ is an indepenent dummy stochastic process with
the same law as $\bm{X}_{t}$. But then 
\[
q(\bm{x})=\P(\bm{X}_{\Delta t}=c_{B})+\P(\bm{X}_{\Delta t}\notin\{c_{A},c_{B}\})\,\E\left[q(\bm{X}_{\Delta t})\right]+O(\Delta t^{2}).
\]
 Expanding further we obtain 
\begin{eqnarray}
 &  & q(\bm{x})=g(\bm{x})\,\Delta t \label{eq:qExpand} \\
 &  & \quad\quad+\ \left(1-f(\bm{x})\,\Delta t-g(\bm{x})\,\Delta t\right)\,\E\left[q\left(\bm{x}-\nabla V(\bm{x})\,\Delta t+\sqrt{2\beta^{-1}\Delta t}\bm{Z}_{0}\right)\right]+O(\Delta t^{2}) \nonumber
\end{eqnarray}
 Then we can expand $q$ via Taylor expansion: 
\begin{eqnarray*}
 &  & q\left(\bm{x}-\nabla V(\bm{x})\,\Delta t+\sqrt{2\beta^{-1}\Delta t}\bm{Z}_{0}\right) = q(\bm{x})-\nabla V(\bm{x})\cdot\nabla q(\bm{x})\,\Delta t \\
 &  & \quad \quad \quad \quad + \ \sqrt{2\beta^{-1}\Delta t}\,\bm{Z}_{0}\cdot\nabla V(\bm{x})+\beta^{-1}\Delta t\,\bm{Z}_{0}^{\top}\nabla^{2}q(\bm{x})\bm{Z}_{0}+o(\Delta t),
\end{eqnarray*}
 from which we obtain 
\begin{eqnarray*}
 &  & \E\left[q\left(\bm{x}-\nabla V(\bm{x})\,\Delta t+\sqrt{2\beta^{-1}\Delta t}\bm{Z}_{0}\right)\right]\\
 &  & \quad \quad=\ \ q(\bm{x})+\left(-\nabla V(\bm{x})\cdot\nabla q(\bm{x})+\beta^{-1}\Delta q(\bm{x})\right)\Delta t+o(\Delta t).
\end{eqnarray*}
It follows from plugging into (\ref{eq:qExpand}) that 
\[
q(\bm{x})=q(\bm{x})+\left[g(\bm{x})-\nabla V(\bm{x})\cdot\nabla q(\bm{x})+\beta^{-1}\Delta q(\bm{x})-(f(\bm{x})+g(\bm{x}))q(\bm{x})\right]\Delta t+o(\Delta t).
\]
Cancelling $q(\bm{x})$ from boths sides, dividing by $\Delta t$,
and taking the limit as $\Delta t\ra0$, we obtain precisely (\ref{eq:softEL}),
as desired.

\section{Discretization of the Ginzburg-Landau density}\label{app:GL-density} In this section, we show that one can approximate the equilibrium distribution of the Ginzburg-Landau potential as \eqref{eq:p-FTT} via the eigenfunctions of a certain kernel.

The computation of the committor function for the Ginzburg-Landau potential requires the numerical approximation of the operator $\mathcal{H}: L^2([-R,R])^d \to \mathbb{R}$ defined by
\begin{align*}
\mathcal{H}[\phi_1,\dots,\phi_d] =& c_\lambda \int_{-R}^R \dots \int_{-R}^R K(0,x_1)\,K(x_1,x_2)\,K(x_2,x_3)\,\dots \, K(x_{d-1},x_d)\\
&\,\times K(x_d,0) \,\phi_1(x_1)\dots\phi_d(x_d)\,{ d}x_1\,\dots\,{ d}x_d,
\end{align*}
where
$$K(x,y) = e^{-\frac{1}{8\lambda}(1-x^2)^2}e^{\frac{\lambda}{2h^2} (x-y)^2}e^{-\frac{1}{8\lambda}(1-y^2)^2},$$
 $c_\lambda = e^{-\frac{1}{4\lambda}}$, and $R$ is some fixed positive constant.

Considering the operator $\mathcal{H}_0 : L^2([-R,R]) \to L^2([-R,R])$ defined by
$$\mathcal{H}_0[\phi](x)=\int_\Omega K(x,y) \,\phi(y)\,{ d}y,$$
we observe that it is compact, symmetric, and positive semi-definite. In particular, it has an eigendecomposition consisting of a countable basis of orthonormal eigenfunctions $u_1,u_2,\dots \in L^2([-R,R])$, together with corresponding non-negative eigenvalues $\lambda_1 \ge \lambda_2 \ge \dots \ge \lambda_n \ge \dots \ge 0$, such that
$$\int_{-R}^R K(x,y)\,\phi(y)\,{ d}y = \sum_{j=1}^\infty \lambda_j u_j(x) \int_{-R}^R u_j(y)\,\phi(x)\,{ d}y$$
for all $\phi \in L^2([-R,R])$. Moreover, for the kernel $K$ defined above it is easily shown that $\lambda_j \in o(e^{-\alpha j})$ for any $\alpha >0$ as $j \to \infty$ (see \cite{LittleReade1984} for example). We note that the implicit constant in the previous estimate will depend on $\alpha$ but not on $j$, and increases rapidly as $\alpha \to \infty$. For convenience, let us define the re-scaled eigenfunctions $v_i = \sqrt{\lambda_i}u_i.$ Upon substitution of the eigendecomposition of $\mathcal{H}_0$ into the definition of $\mathcal{H},$ we obtain
\begin{align}
\label{eqn:H_with_lam}
\mathcal{H}[\phi_1,\dots,\phi_d] =c_\lambda \sum_{j_0,\dots,j_{d}=1}^\infty v_{j_0}(0) A_{j_0,j_1}[\phi_1] \cdots A_{j_{d-1},j_{d}}[\phi_d] v_{j_d}(0),
\end{align}
where $A_{j,\ell}[\phi] = \int_\Omega v_j(x)v_\ell(x) \phi(x)\,{ d}x.$ We observe that $$|A_{j,\ell}[\phi]| \le \sqrt{\lambda_j\lambda_\ell} \|\phi\|_{L^\infty},$$
and thus, if the sums in (\ref{eqn:H_with_lam}) are truncated at term $J$, then
\begin{align*}
&\left|\mathcal{H}[\phi_1,\dots,\phi_d]-\mathcal{H}^{(J)}[\phi_1,\dots,\phi_d]\right|\le \\
&\quad 2\, c_\lambda  \|\phi_1\|_{L^\infty}\cdots \|\phi_d\|_{L^\infty} \sqrt{K(0,0) K^{(j)}(0,0)} { Tr}(\mathcal{H}_0)^{d-\frac{1}{2}} \sqrt{\sum_{j=J+1}^\infty\lambda_j}+\\
&\quad \quad (d-1)\, c_\lambda  \|\phi_1\|_{L^\infty}\cdots \|\phi_d\|_{L^\infty} K(0,0) {\rm Tr}(\mathcal{H}_0)^{d-1} \sum_{j=J+1}^\infty\lambda_j,
\end{align*}
where $K^{(J)}(0,0) := \sum_{j=J+1}^\infty v_j(0)^2 \le K(0,0)$ and $\mathcal{H}^{(J)}$ denotes the truncation of $\mathcal{H}.$ The analyticity of $K$ guarantees that $\sum_{j=J+1}^\infty \lambda_j$ go to zero exponentially quickly \cite{LittleReade1984}) and hence for any $\alpha>0$ there exists a constant $M_{\alpha,\lambda, d,R}$ depending on  $\alpha,$ $\lambda,$ $d,$ and $R$ such that
$$\left|\mathcal{H}[\phi_1,\dots,\phi_d]-\mathcal{H}^{(J)}[\phi_1,\dots,\phi_d]\right| \le M_{\lambda,d,R} \,e^{-\alpha J/2} \|\phi_1\|_{L^\infty} \cdots \|\phi_d\|_{L^\infty}$$
for all $J \ge 1.$

Next, we express the input functions $\phi_i$ in a basis of (suitably-scaled) Chebyshev polynomials,
$$\phi_i(x) = \sum_{n=0}^\infty \phi_{i,n}T_n\left(\frac{x}{R}\right),$$
where $T_n$ is the $n$th standard Chebyshev polynomial. Let $\phi_{i}^{(N)}$ denote the $N$th order truncation of $\phi_i$ defined by
$$\phi_i^{(N)}(x) := \sum_{n=0}^N \phi_{i,n} T_n\left(\frac{x}{R}\right).$$

Then, substituting these Chebyshev expansions into our expression for $\mathcal{H},$ we find
\begin{align}\label{eqn:H_with_cheb}
\mathcal{H}[\phi_1,\dots,\phi_d] =c_\lambda  \sum_{j_0,\dots,j_d =1}^\infty \sum_{n_1,\dots,n_d=1}^\infty A_{j_0,j_1}^{n_1} \cdots A_{j_{d-1},j_d}^{n_d} v_{j_0}(0)v_{j_d}(0)\, \phi_{1,n_1}\dots \phi_{d,n_d},
\end{align}
where
$$A_{j,\ell}^n = A_{j,\ell}\left[ T_n(x/R)\right] .$$
In the following, we assume that $\phi_1,\dots,\phi_d$ are in $C^{p+1}([-R,R])$ for some fixed integer $p \ge 1$ and set 
$$V_p = \max_{i} \left\|  \frac{{ d}^{p+1}}{{ d}x^{p+1}} \phi_i(x/R) \right\|_{L^1([-1,1])}.$$
A standard estimate from approximation theory \cite{Trefethen2013} gives the following bound on the rate of decay of the coefficients of $\phi_1,\dots,\phi_d$:
$$|\phi_{i,n}| \le \frac{2V_p}{\pi n (n-1) \cdots (n-p)}. $$
In particular,
$$\sum_{n=N+1}^\infty |\phi_{i,n}| \le \frac{2V_p}{\pi p (N-p)^p}.$$

If the sums over the Chebyshev coefficients (the $n$ indices)  in (\ref{eqn:H_with_cheb}) are truncated at a fixed integer $N$ and the sums over the eigenvalues (the $j$ indices) are truncated at $J$, then the error is bounded by
\begin{align*}
&\left| \mathcal{H}[\phi_1,\dots,\phi_d] - \mathcal{H}^{(J,N)}[\phi_1,\dots,\phi_d]\right|\le M_{\lambda,d,R}\, e^{-\alpha J/2} \|\phi_1\|_{L^\infty} \cdots \|\phi_d\|_{L^\infty} \\
&\quad + \frac{2 d c_\lambda V_p}{\pi p (N-p)^p}  \left(1+\frac{2V_p}{\pi p (N-p)^p} \right)^{d-1} K(0,0)\, {\rm Tr}(\mathcal{H}_0)^{d} \max_i \|\phi_i\|_{L^\infty}^{d-1}
\end{align*}
Here $\mathcal{H}^{(J,N)}$ denotes the operator obtained by truncating the sums over eigenvalues in $\mathcal{H}$ at $J$ and projecting onto the first $N+1$ terms in the Chebyshev expansions of $\phi_1,\dots,\phi_d.$ This latter projection, along with the eigendecomposition of $\mathcal{H}_0,$ can be performed easily on the computer using standard numerical integration. This yields a discrete, finite-dimensional tensor which is the object we use in our approach when approximating $\mathcal{H}.$

\newpage
\bibliographystyle{siam}
\bibliography{ref}

\begin{thebibliography}{10}

\bibitem{mps-1}
{\sc Ian Affleck, Tom Kennedy, Elliott~H Lieb, and Hal Tasaki}, {\em Valence
  bond ground states in isotropic quantum antiferromagnets}, in Condensed
  matter physics and exactly soluble models, Springer, 1988, pp.~253--304.

\bibitem{tensor-pde-2}
{\sc Markus Bachmayr, Reinhold Schneider, and Andr{\'e} Uschmajew}, {\em Tensor
  networks and hierarchical tensors for the solution of high-dimensional
  partial differential equations}, Foundations of Computational Mathematics, 16
  (2016), pp.~1423--1472.

\bibitem{conf-change-1}
{\sc Anna Berteotti, Andrea Cavalli, Davide Branduardi, Francesco~Luigi
  Gervasio, Maurizio Recanatini, and Michele Parrinello}, {\em Protein
  conformational transitions: the closure mechanism of a kinase explored by
  atomistic simulations}, Journal of the American Chemical Society, 131 (2009),
  pp.~244--250.

\bibitem{bhattacharya2009stochastic}
{\sc Rabi~N Bhattacharya and Edward~C Waymire}, {\em Stochastic processes with
  applications}, SIAM, 2009.

\bibitem{foundation-of-ds}
{\sc Avrim Blum, John Hopcroft, and Ravindran Kannan}, {\em Foundations of data
  science}, Vorabversion eines Lehrbuchs, 5 (2016).

\bibitem{diffusion-map-committor}
{\sc Ronald~R Coifman, Ioannis~G Kevrekidis, St{\'e}phane Lafon, Mauro
  Maggioni, and Boaz Nadler}, {\em Diffusion maps, reduction coordinates, and
  low dimensional representation of stochastic systems}, Multiscale Modeling \&
  Simulation, 7 (2008), pp.~842--864.

\bibitem{tt-quantum-1}
{\sc David~Elieser Deutsch}, {\em Quantum computational networks}, Proceedings
  of the Royal Society of London. A. Mathematical and Physical Sciences, 425
  (1989), pp.~73--90.

\bibitem{tt-quantum-5}
{\sc Glen Evenbly and Guifr{\'e} Vidal}, {\em Tensor network states and
  geometry}, Journal of Statistical Physics, 145 (2011), pp.~891--918.

\bibitem{conf-change-2}
{\sc Barry~J Grant, Alemayehu~A Gorfe, and J~Andrew McCammon}, {\em Large
  conformational changes in proteins: signaling and other functions}, Current
  opinion in structural biology, 20 (2010), pp.~142--147.

\bibitem{tt-review-3}
{\sc Lars Grasedyck}, {\em Hierarchical low rank approximation of tensors and
  multivariate functions}, Lecture notes of the Z{\"u}rich summer school on
  Sparse Tensor Discretizations of High-Dimensional Problems,  (2010).

\bibitem{tt-review-2}
{\sc Lars Grasedyck, Daniel Kressner, and Christine Tobler}, {\em A literature
  survey of low-rank tensor approximation techniques}, GAMM-Mitteilungen, 36
  (2013), pp.~53--78.

\bibitem{h-tucker}
{\sc Wolfgang Hackbusch, Boris~N Khoromskij, and Eugene~E Tyrtyshnikov}, {\em
  Hierarchical kronecker tensor-product approximations},  (2005).

\bibitem{deep-pde-2}
{\sc Jiequn Han, Arnulf Jentzen, and E~Weinan}, {\em Solving high-dimensional
  partial differential equations using deep learning}, Proceedings of the
  National Academy of Sciences, 115 (2018), pp.~8505--8510.

\bibitem{cp-decomposition}
{\sc Frank~L Hitchcock}, {\em The expression of a tensor or a polyadic as a sum
  of products}, Journal of Mathematics and Physics, 6 (1927), pp.~164--189.

\bibitem{GL-model}
{\sc K-H Hoffmann and Qi~Tang}, {\em Ginzburg-Landau phase transition theory
  and superconductivity}, vol.~134, Birkh{\"a}user, 2012.

\bibitem{NN-committor-1}
{\sc Yuehaw Khoo, Jianfeng Lu, and Lexing Ying}, {\em Solving for
  high-dimensional committor functions using artificial neural networks},
  Research in the Mathematical Sciences, 6 (2019), pp.~1--13.

\bibitem{deep-pde-3}
\leavevmode\vrule height 2pt depth -1.6pt width 23pt, {\em Solving parametric
  pde problems with artificial neural networks}, European Journal of Applied
  Mathematics, 32 (2021), pp.~421--435.

\bibitem{tt-review-4}
{\sc Boris~N Khoromskij}, {\em Tensor numerical methods for multidimensional
  pdes: theoretical analysis and initial applications}, ESAIM: Proceedings and
  Surveys, 48 (2015), pp.~1--28.

\bibitem{tensorGalerkin2011}
{\sc Boris~N. Khoromskij and Christoph Schwab}, {\em Tensor-structured galerkin
  approximation of parametric and stochastic elliptic pdes}, SIAM Journal on
  Scientific Computing, 33 (2011), pp.~364--385.

\bibitem{tensor-pde-1}
{\sc Boris~N Khoromskij and Christoph Schwab}, {\em Tensor-structured galerkin
  approximation of parametric and stochastic elliptic pdes}, SIAM Journal on
  Scientific Computing, 33 (2011), pp.~364--385.

\bibitem{tt-review-1}
{\sc Tamara~G Kolda and Brett~W Bader}, {\em Tensor decompositions and
  applications}, SIAM review, 51 (2009), pp.~455--500.

\bibitem{point-committor}
{\sc Rongjie Lai and Jianfeng Lu}, {\em Point cloud discretization of
  fokker--planck operators for committor functions}, Multiscale Modeling \&
  Simulation, 16 (2018), pp.~710--726.

\bibitem{chem-tran-1}
{\sc Antonio~C Lasaga}, {\em 2. transition state theory}, in Kinetic Theory in
  the Earth Sciences, Princeton University Press, 2014, pp.~152--219.

\bibitem{NN-committor-2}
{\sc Haoya Li, Yuehaw Khoo, Yinuo Ren, and Lexing Ying}, {\em Solving for high
  dimensional committor functions using neural network with online
  approximation to derivatives}, arXiv preprint arXiv:2012.06727,  (2020).

\bibitem{LiLinRen2019}
{\sc Qianxiao Li, Bo~Lin, and Weiqing Ren}, {\em Computing committor functions
  for the study of rare events using deep learning}, The Journal of Chemical
  Physics, 151 (2019), p.~054112.

\bibitem{LittleReade1984}
{\sc G.~Little and J.B. Reade}, {\em Eigenvalues of analytic kernels}, SIAM J.
  Math. Anal., 15 (1984), pp.~133--136.

\bibitem{transition-2}
{\sc Jianfeng Lu and James Nolen}, {\em Reactive trajectories and the
  transition path process}, Probability Theory and Related Fields, 161 (2015),
  pp.~195--244.

\bibitem{chem-tran-2}
{\sc Naoto Okuyama-Yoshida, Masataka Nagaoka, and Tokio Yamabe}, {\em
  Transition-state optimization on free energy surface: Toward solution
  chemical reaction ergodography}, International journal of quantum chemistry,
  70 (1998), pp.~95--103.

\bibitem{tt-quantum-4}
{\sc Rom{\'a}n Or{\'u}s}, {\em Advances on tensor network theory: symmetries,
  fermions, entanglement, and holography}, The European Physical Journal B, 87
  (2014), pp.~1--18.

\bibitem{tensor-diagram-intro}
\leavevmode\vrule height 2pt depth -1.6pt width 23pt, {\em A practical
  introduction to tensor networks: Matrix product states and projected
  entangled pair states}, Annals of Physics, 349 (2014), pp.~117--158.

\bibitem{tt-decomposition}
{\sc Ivan~V Oseledets}, {\em Tensor-train decomposition}, SIAM Journal on
  Scientific Computing, 33 (2011), pp.~2295--2317.

\bibitem{mps-2}
{\sc David Perez-Garcia, Frank Verstraete, Michael~M Wolf, and J~Ignacio
  Cirac}, {\em Matrix product state representations}, arXiv preprint
  quant-ph/0608197,  (2006).

\bibitem{pinn}
{\sc Maziar Raissi, Paris Perdikaris, and George~E Karniadakis}, {\em
  Physics-informed neural networks: A deep learning framework for solving
  forward and inverse problems involving nonlinear partial differential
  equations}, Journal of Computational Physics, 378 (2019), pp.~686--707.

\bibitem{transition-info-2}
{\sc Weiqing Ren, Eric Vanden-Eijnden, Paul Maragakis, and Weinan E}, {\em
  Transition pathways in complex systems: Application of the finite-temperature
  string method to the alanine dipeptide}, The Journal of chemical physics, 123
  (2005), p.~134109.

\bibitem{Rostkoff2021}
{\sc Grant~M. Rostkoff, Andrew~R. Mitchell, and Eric Vanden-Eijnden}, {\em
  Active importance sampling for variational objectives dominated by rare
  events: Consequences for optimization and generalization}, arXiv:2008.06334.

\bibitem{deep-pde-1}
{\sc Justin Sirignano and Konstantinos Spiliopoulos}, {\em Dgm: A deep learning
  algorithm for solving partial differential equations}, Journal of
  computational physics, 375 (2018), pp.~1339--1364.

\bibitem{dga}
{\sc Erik~H. Thiede, Dimitrios Giannakis, Aaron~R. Dinner, and Jonathan Weare},
  {\em Galerkin approximation of dynamical quantities using trajectory data},
  The Journal of Chemical Physics, 150 (2019), p.~244111.

\bibitem{Trefethen2013}
{\sc L.~N. Trefethen}, {\em Approximation {T}heory and {A}pproximation
  {P}ractice}, SIAM, 2013.

\bibitem{tucker-decomposition}
{\sc Ledyard~R Tucker}, {\em Some mathematical notes on three-mode factor
  analysis}, Psychometrika, 31 (1966), pp.~279--311.

\bibitem{transition-3}
{\sc Eric Vanden-Eijnden et~al.}, {\em Transition-path theory and path-finding
  algorithms for the study of rare events.}, Annual review of physical
  chemistry, 61 (2010), pp.~391--420.

\bibitem{revisit-finite-temp-string}
{\sc Eric Vanden-Eijnden and Maddalena Venturoli}, {\em Revisiting the finite
  temperature string method for the calculation of reaction tubes and free
  energies}, The Journal of chemical physics, 130 (2009), p.~05B605.

\bibitem{vershynin2018high}
{\sc Roman Vershynin}, {\em High-dimensional probability: An introduction with
  applications in data science}, vol.~47, Cambridge university press, 2018.

\bibitem{tt-quantum-3}
{\sc Frank Verstraete, Valentin Murg, and J~Ignacio Cirac}, {\em Matrix product
  states, projected entangled pair states, and variational renormalization
  group methods for quantum spin systems}, Advances in Physics, 57 (2008),
  pp.~143--224.

\bibitem{transition-info-1}
{\sc E~Weinan, Weiqing Ren, and Eric Vanden-Eijnden}, {\em Transition pathways
  in complex systems: Reaction coordinates, isocommittor surfaces, and
  transition tubes}, Chemical Physics Letters, 413 (2005), pp.~242--247.

\bibitem{finite-temp-string}
{\sc E~Weinan, Weiqing Ren, Eric Vanden-Eijnden, et~al.}, {\em Finite
  temperature string method for the study of rare events}, J. Phys. Chem. B,
  109 (2005), pp.~6688--6693.

\bibitem{transition-1}
{\sc E~Weinan and Eric Vanden-Eijnden}, {\em Towards a theory of transition
  paths}, Journal of statistical physics, 123 (2006), pp.~503--523.

\bibitem{deepritz}
{\sc E~Weinan and Bing Yu}, {\em The deep ritz method: a deep learning-based
  numerical algorithm for solving variational problems}, Communications in
  Mathematics and Statistics, 6 (2018), pp.~1--12.

\bibitem{mps-3}
{\sc Steven~R White}, {\em Density matrix formulation for quantum
  renormalization groups}, Physical review letters, 69 (1992), p.~2863.

\bibitem{nuc-1}
{\sc Dirk Zahn and Stefano Leoni}, {\em Nucleation and growth in
  pressure-induced phase transitions from molecular dynamics simulations:
  Mechanism of the reconstructive transformation of nacl to the cscl-type
  structure}, Physical review letters, 92 (2004), p.~250201.

\bibitem{nuc-2}
{\sc Bingge Zhao, Linfang Li, Fenggui Lu, Qijie Zhai, Bin Yang, Christoph
  Schick, and Yulai Gao}, {\em Phase transitions and nucleation mechanisms in
  metals studied by nanocalorimetry: A review}, Thermochimica Acta, 603 (2015),
  pp.~2--23.

\end{thebibliography}

\end{document}